\documentclass{amsart}
\title{Combinatorial $N_\infty$ operads}

\author{Jonathan Rubin}
\address{University of California Los Angeles,
Los Angeles, CA 90095}
\email{jrubin@math.ucla.edu}


\date{\today}

\usepackage{amsmath,amssymb,amsthm,stackrel,tikz,enumitem,stmaryrd,mathrsfs,physics}

\newcommand{\s}[1]{\mathscr{#1}}

\renewcommand{\t}[1]{\textnormal{#1}}
\newcommand{\ol}[1]{\overline{#1}}
\newcommand{\ul}[1]{\underline{#1}}
\renewcommand{\c}[1]{\mathcal{#1}}
\renewcommand{\b}[1]{\mathbf{#1}}
\newcommand{\ub}[1]{\underline{\mathbf{#1}}}
\newcommand{\til}[1]{\widetilde{#1}}
\newcommand{\bb}[1]{\mathbb{#1}}

\newcommand{\ve}{\varepsilon}
\newcommand{\vp}{\varphi}
\newcommand{\la}{\langle}
\newcommand{\ra}{\rangle}

\theoremstyle{definition}
\newtheorem{defn}{Definition}[section]
\newtheorem{ex}[defn]{Example}
\newtheorem{cex}[defn]{Counterexample}

\newtheorem{const}[defn]{Construction}
\newtheorem{conj}[defn]{Conjecture}
\newtheorem{nota}[defn]{Notation}

\theoremstyle{plain}
\newtheorem{lem}[defn]{Lemma}
\newtheorem{thm}[defn]{Theorem}
\newtheorem{prop}[defn]{Proposition}
\newtheorem{cor}[defn]{Corollary}

\newtheorem*{thm1}{\protect{Theorem \ref{thm:Gsetmod}}}
\newtheorem*{thm2}{\protect{Theorem \ref{thm:freerealize}}}

\newtheorem*{thm4}{\protect{Theorem \ref{thm:modOpG+} and Proposition \ref{prop:LHNOp}}}
\newtheorem*{thm5}{\protect{Theorem \ref{thm:assocNop}}}
\newtheorem*{thmkeycalc*}{\protect{Theorem \ref{thm:admfree}}}
\newtheorem*{conj*}{Conjecture}

\theoremstyle{remark}
\newtheorem{rem}[defn]{Remark}

\begin{document}
\maketitle

\begin{abstract} We prove that the homotopy theory of $N_\infty$ operads is equivalent to a homotopy theory of discrete operads, and we construct free and associative operadic realizations of every indexing system. This resolves a conjecture of Blumberg and Hill in the affirmative.
\end{abstract}

\tableofcontents

\section{Introduction}\label{sec:intro}

Operads were first introduced by May in \cite{MayGILS}, and they have been applied throughout algebra and topology ever since. As the name might suggest, an operad is an object that parametrizes operations. They appear in many contexts, and interesting structure on an operad $\s{O}$ translates universally into interesting structure on the algebras over $\s{O}$.

The original application of operad theory was the recognition principle for iterated loop spaces, due to Boardman-Vogt, May, and Milgram. May's approach to this theorem leverages operadic structure on $X$ to construct an equivalence between $X$ and an $n$-fold loop space $\Omega^n Y = \t{Map}_*(S^n,Y)$. The basic idea is to track the homotopy coherence of the sum $+$ in $\pi_n(Y) = \pi_0(\Omega^n Y)$. Two maps $f , g : S^n = I^n/\partial I^n  \rightrightarrows Y$ are usually added together by pasting $f$ and $g$ onto two different halves of the $n$-cube $I^n$, but there are $n$ different dimensions to choose from, and this operation is only a group structure up to homotopy. Moreover, there are many homotopies that witness the associativity, unitality, and the commutativity of $+$ when $n > 1$. The little $n$-cubes operad $\c{C}_n$ parametrizes all of the possibilities, and the operadic recognition principle states that $X$ is a $n$-fold loop space if and only if $X$ is a grouplike $\c{C}_n$-algebra, i.e. $\pi_0(X)$ is a group.

As $n$ increases, there are more and more degrees of freedom in $I^n$, and the connectivity of the operads $\c{C}_n$ correspondingly increases. Passing to colimits yields the infinite little cubes operad $\c{C}_\infty$, which parametrizes the additive structure on infinite loop spaces. This is the prototypical example of an $E_\infty$ operad. It parametrizes operations that are associative, unital, and commutative, up to all possible homotopies, and its actions can be used to construct infinite deloopings.

Now, grouplike infinite loop spaces are equivalent to connective spectra, with the $\c{C}_\infty$-action corresponding to addition. On the other hand, multiplicative structures on spectra are classically parametrized by a different operad. If $E$ is a spectrum indexed over the subspaces of $\bb{R}^\infty$, then $E^{\land n}$ is naturally indexed over the subspaces of $(\bb{R}^\infty)^{\oplus n}$. After changing universe along a linear isometry $f : (\bb{R}^{\infty})^{\oplus n} \to \bb{R}^\infty$, we can map back to $E$, but there are many possible choices for $f$. The linear isometries operad $\c{L}$ parametrizes all of the options. It is also an $E_\infty$ operad, but its geometry differs greatly from that of $\c{C}_\infty$. Nevertheless, there is a zig-zag of equivalences $\c{C}_\infty \stackrel{\sim}{\leftarrow} \c{C}_\infty \times \c{L} \stackrel{\sim}{\to} \c{L}$ connecting them. This is May's ``product trick.'' It implies that all $E_\infty$ operads are equivalent. 

The situation is not nearly so clear-cut in equivariant homotopy theory. Suppose $G$ is a finite group. Then there are $G$-equivariant analogues to the operads $\c{C}_\infty$ and $\c{L}$. The equivariant version of $\c{L} = \c{L}(\bb{R}^\infty)$ is obtained by replacing $\bb{R}^\infty$ with a $G$-universe $U$. We think of $\c{L}(U)$ as the natural representing object for multiplication on $G$-spectra over $U$. The equivariant version of $\c{C}_\infty$ is more subtle. Cubes are ``too square'' to support a $G$-action, so one replaces the cubes in $\c{C}_n$ with the unit discs of finite-dimensional $G$-representations $V$. The result is the little $V$-discs operad $\c{D}(V)$. Given a universe $U$, one takes a colimit over finite-dimensional subrepresentations $V \subset U$ to get the infinite little discs operad $\c{D}(U) = \t{colim}_{V \subset U} \c{D}(V)$, but this does not naturally act on equivariant iterated loop spaces. However, there is a thickening of $\c{D}(U)$ which does act, namely the Steiner operad $\c{K}(U)$. We think of $\c{K}(U)$ as the natural representing object for addition on $G$-spectra over $U$.

Now suppose that $R$ is a genuine commutative ring $G$-spectrum. Ignoring multiplication for a moment, the $\bb{Z}$-graded homotopy groups of $R$ are naturally $G$-Mackey functors. We may understand their transfers in terms of the $\c{K}(U)$-action. Indeed, additive transfers are usually constructed by embedding an orbit into a representation $V$ and then taking the Pontrjagin-Thom collapse map. This corresponds to an operation in $\c{D}(U)$ and also in $\c{K}(U) \simeq \c{D}(U)$. On the other hand, there are also multiplicative norms in the $RO(G)$-graded homotopy of $R$, first introduced by Greenlees and May \cite{GreMay}, and used to great effect by Hill-Hopkins-Ravenel \cite{HHR}. One can similarly understand these norms in terms of the $\c{L}(U)$-action. On the level of universes, norms from $H$ to $G$ arise from certain $G$-equivariant linear isometries $f : U^{\oplus n} \to U$, for which the $G$-action on $U^{\oplus n}$ is restricted from an action of $\Sigma_n \wr H$.

Thus, equivariant $E_\infty$ operads parametrize much more than just a homotopy coherent commutative monoid operation $*$. They also parametrize transfers or norms, depending on whether we think of $*$ as additive or multiplicative. If $U$ is a complete universe, then $\c{K}(U)$ and $\c{L}(U)$ parametrize all transfers and norms, and if $U$ is a trivial universe, then $\c{K}(U)$ and $\c{L}(U)$ parametrize no transfers or norms. A surprising observation, due to Blumberg and Hill \cite[Theorem 4.22]{BH1}, is that there are incomplete universes $U$ such that $\c{K}(U)$ and $\c{L}(U)$ parametrize different sets of transfers and norms. Thus, as $U$ varies over all possible $G$-universes, we obtain distinct families of operads $\c{K}(U) \simeq \c{D}(U)$ and $\c{L}(U)$. These are the prototypical examples of $N_\infty$ operads.

In general, a $N_\infty$ operad is a $G$-equivariant operad that parametrizes a homotopy coherent commutative monoid structure together with a compatible system of (additive or multiplicative) transfers. When $G$ is the trivial group, a $N_\infty$ operad is just an $E_\infty$ operad, and all $E_\infty$ operads are equivalent. For general groups $G$, the homotopy type of a $N_\infty$ $G$-operad is completely determined by its transfers. As explained above, there are multiple possibilities, so it makes sense to try to classify them. Blumberg and Hill began such a program in \cite{BH1}. Given any $N_\infty$ operad $\s{O}$, they construct an ``indexing system,'' which encodes the transfers of $\s{O}$. This is a combinatorial object, which satisfies axioms that encode how transfers interact with an operad structure. For any group $G$, the collection of all $G$-indexing forms a lattice under inclusion, and maps $\s{O}_1 \to \s{O}_2$ between $N_\infty$ operads induce inclusions of indexing systems. Thus, we obtain a functor from the category $N_\infty\t{-}\b{Op}^G$ of $N_\infty$ $G$-operads to the poset category $\b{Ind}(G)$ of all $G$-indexing systems. In fact, this functor factors through the homotopy category $\t{Ho}(N_\infty\t{-}\b{Op}^G)$ because equivalent $N_\infty$ operads have equal indexing systems.

Blumberg and Hill proved that $\t{Ho}(N_\infty\t{-}\b{Op}^G)$ is mapped fully and faithfully into $\b{Ind}(G)$. They also made the following

\begin{conj*} Taking indexing systems determines an equivalence between the category $\t{Ho}(N_\infty\t{-}\b{Op}^G)$ and the poset $\b{Ind}(G)$.
\end{conj*}

In other words, Blumberg and Hill conjectured that every indexing system is realized by some $N_\infty$ operad. In this paper, we shall give a combinatorial verification Blumberg and Hill's conjecture. Other solutions to this problem have been found independently by Guti\'{e}rrez and White \cite{GutWhite}, and by Bonventre and Pereira \cite{BonPer}, and we give a quick comparison between our constructions in \S\ref{subsec:compare}.

Our three solutions are very different, and they highlight complementary aspects of equivariant operad theory. Guti\'{e}rrez and White study a myriad of model category structures on the category of $G$-operads, much in the spirit of Berger-Moerdijk \cite{BerMoer}. Their realizations of indexing systems arise as cofibrant replacements of the commutativity operad in judiciously chosen model categories. In contrast, Bonventre and Pereira introduce a novel kind of equivariant operad, which are a blend of ordinary operads and fixed-point presheaves. Thus, they build norms into the underlying formalism, and their realizations of indexing systems arise as operadic variants of Elmendorf's construction of universal spaces \cite{Elm}.

The purpose of this paper is to reduce $N_\infty$ theory to combinatorics. This drastically simplifies the mathematics, and it brings precise, algebraic theorems within arm's reach. In effect, our approach strips away all of the topology, leaving only the algebra of discrete equivariant operads. That being said, this algebra is rather nontrivial. An operad is a generalization of a monoid, and the most interesting operads arise as quotients. Thus, we are forced to contend with word problems. Of course, these word problems are also present in the topological case, but our work demonstrates that they are, in some sense, the only problems.

More precisely, we introduce discrete analogues to $N_\infty$ operads, which we call $N$ operads, and then we prove the following result.

\begin{thm1}The category of $N_\infty$ operads and the category of $N$ operads have equivalent hammock localizations, and this equivalence respects indexing systems.
\end{thm1}

This theorem can be refined. The category of $N_\infty$ operads is not bicomplete, and therefore cannot admit a model category structure, but we can replace $N$ operads with a model category that has the same underlying homotopy theory. Let $\b{Op}^G_+$ denote the category of operads in $G$-sets, equipped with a marked $G$-fixed constant and $G$-equivariant binary product.

\begin{thm4} The category $\b{Op}^G_+$ supports a right proper, combinatorial, simplicial model category structure. This model category has the same hammock localization as the category of $N_\infty$ operads.
\end{thm4}

This model category structure on $\b{Op}^G_+$ has a number of uses. Looking ahead, it is indispensible in \cite{RubCateg}, where we lift natural operations on indexing systems back to the operad level. In this paper, we use it to give a new, combinatorial proof that $\t{Ho}(N_\infty\t{-}\b{Op}^G)$ embeds into $\b{Ind}(G)$, and we also use it to contexualize our first major construction.

\begin{thm2} Every indexing system $\c{I}$ is realizable by a finitely generated free $N$ operad $\b{F}(\c{I})$, which may be constructed functorially in $\c{I}$.
\end{thm2}

From a conceptual standpoint, the operad $\b{F}(\c{I})$ is a cofibrant replacement of the commutativity operad in a suitable model structure on $\b{Op}^G_+$ (cf. Proposition \ref{prop:freerealcof}). This is formally analogous to the situation in \cite{GutWhite} and \cite{BonPer}, and after passing to $N_\infty$ operads, we obtain a similar operad to theirs (cf. \S\ref{subsec:compare}). Theorem \ref{thm:freerealize} resolves Blumberg and Hill's conjecture, but it also goes further.

For example, the finite generation of $\b{F}(\c{I})$ is of great use. One can construct a categorical $N_\infty$ operad $\til{\b{F}}(\c{I})$ by applying the right adjoint to the object functor $\t{Ob} : \b{Cat} \to \b{Set}$, and we prove in \cite{RubNSMC} that $\til{\b{F}}(\c{I})$-algebra $G$-categories are ``normed symmetric monoidal categories'' (NSMCs), i.e. ordinary symmetric monoidal categories equipped with certain twisted products. The finite generation of $\b{F}(\c{I})$ ensures that NSMCs are finitely presentable, which is in sharp contrast to Guillou-May-Merling-Osorno's symmetric monoidal $G$-categories (cf. \cite{GMMO} and \cite{BBKOOTX}). We emphasize that the finite generation of $\b{F}(\c{I})$ is a consequence of the combinatorics of indexing systems, rather than the model-categorical formalism.

Just as Theorem \ref{thm:Gsetmod} can be refined, so too can Theorem \ref{thm:freerealize}.

\begin{thm5} Every indexing system $\c{I}$ is realizable by a finitely presented, associative and unital $N$ operad $\b{As}(\c{I})$, which may be constructed functorially in $\c{I}$.
\end{thm5}

In contrast to $\b{F}(\c{I})$, the operad $\b{As}(\c{I})$ is invisible to the model-category theory because it is not cofibrant. However, it has a number of convenient properties. To start, it is very small. It has no nontrivial nullary or unary operations, and $\b{As}(\c{I})(n)$ grows far more slowly than the $G$-permutativity operad considered by \cite{GMMO}. Applying the right adjoint to the object functor $\t{Ob} : \b{Cat} \to \b{Set}$ yields a $N_\infty$ permutativity operad $\s{P}(\c{I}) = \til{\b{As}}(\c{I})$, whose algebras are strictly associative and unital NSMCs, i.e. normed permutative categories. We suspect that these structures will be useful in categorical infinite loop space theory, but that remains to be seen. On the other hand, if we pass to space-level $N_\infty$ operads, then we obtain an equivariant Barratt-Eccles operad $\s{E}(\c{I})$. The operad $\s{E}(\c{I})$ is reduced, which is technically convenient in \cite[Remark 2.7]{BH3}. We do not know of any other general construction of reduced $N_\infty$ operads.

\subsection*{Organization} The remainder of this paper is organized as follows. In \S\ref{sec:class}, we give a quick introduction to the theory of $N_\infty$ operads. We recall some basic definitions and examples, and then we summarize the classification theorem. In \S\ref{sec:Nop}, we introduce $N$ operads, explain their relationship to $N_\infty$ operads, and give examples. We prove that $N$ operads and $N_\infty$ operads have the same homotopy theory (Theorem \ref{thm:Gsetmod}). In \S\ref{sec:realize}, we explain how to construct free realizations of indexing systems (Theorem \ref{thm:freerealize}), modulo the calculation of the fixed points of a free operad (Theorem \ref{thm:admfree}). Theorem \ref{thm:admfree} is the key technical result of this paper. We set up some scaffolding in \S\ref{sec:quotfreeop}, and then we do the calculation in \S\ref{subsec:pfkeycalc}. In \S\ref{sec:assocNop}, we introduce associative $N$ operads and establish their basic properties. This strengthens the result in \S\ref{sec:realize}. Lastly, we spend \S\ref{sec:modstr} developing the model category theory of discrete operads in $G$-sets.

The reader who wants a quick introduction to $N_\infty$ theory should read \S\ref{sec:class}. The reader who wants a summary of our solution to Blumberg and Hill's conjecture should read \S\ref{subsec:Nop}, \S\ref{sec:realize}, \S\ref{subsec:sumfree}, and \S\ref{subsec:pfkeycalc}.

\subsection*{Conventions} Throughout this paper, $G$ denotes a finite, discrete group with unit $e$, and all spaces are understood to be compactly generated and weak Hausdorff. All of our operads are symmetric operads in an ambient cartesian monoidal category. Typically, this will be the category of left $G$-spaces or left $G$-sets.

\subsection*{Acknowledgements} It is a pleasure to thank Peter May and Mike Hill for guidance and inspiration throughout this project. We also thank Ang\'{e}lica Osorno, Kyle Ormsby, and the students at Reed for stimulating conversations that prompted this revision. This work was partially supported by NSF Grant DMS-1803426.

\section{The classification of $N_\infty$ operads}\label{sec:class}

This section is a brief introduction to the theory of $N_\infty$ operads. We review some key concepts and examples, and then we summarize the classification of $N_\infty$ operads (Theorem \ref{thm:classNinfty}). Our discussion is based heavily on \cite{BH1} and \cite{GM}. With the exception of the surjectivity portion of Theorem \ref{thm:classNinfty}, the contents of this section were already known. The surjectivity follows independently from \cite{BonPer}, \cite{GutWhite}, and Theorems \ref{thm:freerealize} or \ref{thm:assocNop} in this paper.

\subsection{Equivariant operads}

Let $G$ be a finite group with unit $e$. Throughout this discussion, we work in the category $\b{Top}^G$ of left $G$-spaces and $G$-equivariant continuous maps. The category $\b{Top}^G$ carries two natural enrichments. Let $\ub{To}\b{p}(X,Y)$ denote the space of all continuous maps from $X$ to $Y$ equipped with the compact-open topology. On the one hand, we can topologize the set $\b{Top}^G(X,Y)$ of $G$-equivariant continuous maps $X \to Y$ as a subspace $\ub{To}\b{p}^G(X,Y) \subset \ub{To}\b{p}(X,Y)$, which enriches $\b{Top}^G$ over $\b{Top}$. On the other hand, $\b{Top}^G$ is a cartesian closed category, whose products $X \times Y$ are equipped the diagonal $G$-action, and whose internal homs $\ub{To}\b{p}_G(X,Y)$ are the spaces $\ub{To}\b{p}(X,Y)$ equipped with the conjugation $G$-action. These two enrichments are related through $\ub{To}\b{p}_G(X,Y)^G = \ub{To}\b{p}^G(X,Y)$, but the hom $G$-spaces $\ub{To}\b{p}_G(X,Y)$ are more relevant to operad theory.

The prototypical example of a $G$-operad is the endomorphism operad $\b{End}(X)$ of a $G$-space $X$. The $n$th level of $\b{End}(X)$ is the hom $G$-space $\ub{To}\b{p}_G(X^{\times n},X)$. It carries a left conjugation $G$-action and a right permutation $\Sigma_n$-action, and these actions commute. We usually repackage this structure into a single left $G \times \Sigma_n$-action $(g,\sigma) \cdot f = gf\sigma^{-1}$. The identity map $\t{id} \in \ub{To}\b{p}_G(X,X)$ is $G$-fixed, the composition operation $\gamma(h;f_1,\dots,f_{\t{arity}(h)}) = h \circ (f_1 \times \cdots \times f_{\t{arity}(h)})$ is $G$-equivariant with respect to conjugation, and evident associativity, unitality, and $\Sigma$-equivariance relations hold. This structure is axiomatized in the following definition.

\begin{defn} A \emph{$G$-operad} $\s{O}$ is a symmetric operad in the category $\b{Top}^G$. Explicitly, $\s{O}$ consists of a sequence $(\s{O}(n))_{n \geq 0}$ of $G \times \Sigma_n$-spaces, equipped with a $G$-fixed identity $\t{id} \in \s{O}(1)$ and a continuous $G$-equivariant composition map 
	\[
	\gamma : \s{O}(k) \times \s{O}(j_1) \times \cdots \times \s{O}(j_k) \to \s{O}(j_1 + \cdots + j_k)
	\]
for every $k,j_1,\dots,j_k \geq 0$, such that the usual associativity, unitality, and $\Sigma$-equivariance axioms hold (cf. \cite[Definition 1.1]{MayGILS}). We write $\abs{f}$ for the \emph{arity} of an operation in $\s{O}$. Thus $\abs{f} = n$ means $f \in \s{O}(n)$.

A map $\vp : \s{O}_1 \to \s{O}_2$ of $G$-operads is a sequence of continuous, $G \times \Sigma_n$-equivariant maps $\vp_n : \s{O}_1(n) \to \s{O}_2(n)$ that preserve the identity and composition. An \emph{$\s{O}$-algebra $G$-space} is a representation of $\s{O}$ over a $G$-space, i.e. an object $X \in \b{Top}^G$ equipped with an operad map $\s{O} \to \b{End}(X)$.
\end{defn}

We think of the $n$th level $\s{O}(n)$ of a $G$-operad as a parameter space for $n$-ary operations on a $G$-space $X$. The stabilizer of $f \in \s{O}(n)$ encodes the $G$-equivariance and commutativity relations that $f : X^{\times n} \to X$ satisfies. For example, if $f$ is $G$-fixed, then $f : X^{\times n} \to X$ is $G$-equivariant, and if $f$ is $\Sigma_n$-fixed, then $f(x_{\sigma^{-1}1},\dots,x_{\sigma^{-1}n}) = f(x_1,\dots,x_n)$ for every permutation of its arguments. More interesting relations appear when $G$ and $\Sigma_n$ isotropy conditions mix.

We start with the simplest case. Regard the commutativity operad $\b{Com}$ as a discrete $G$-operad with trivial $G$-action. The $n$th level of $\b{Com}$ is $*$ for all $n \geq 0$. A $\b{Com}$-algebra $G$-space $X$ is a strictly associative, commutative, and unital monoid in $\b{Top}^G$, whose product $*$ is strictly $G$-equivariant, i.e.
	\[
	g(x * y) = (gx) * (gy),
	\]
and whose unit element $1 \in X$ is strictly $G$-fixed. 

Now consider the fixed point subspaces $X^H$ of $X$. Every inclusion $K \subset H$ of subgroups gives a reverse inclusion $X^K \supset X^H$ on fixed points, every element $g \in G$ gives an isomorphism $g \cdot (-) : X^H \to X^{gHg^{-1}}$, and these data determine the equivariant homotopy type of $X$ by Elmendorf's theorem \cite[Theorem 1]{Elm}. However, there is additional structure on the fixed points of $X$ coming from the operad action. Since $*$ is $G$-equivariant and $1$ is $G$-fixed, the monoid structure on $X$ restricts to every subspace $X^H$. More interestingly, for every inclusion $K \subset H$ of subgroups, there is a ``wrong-way'' norm map $\t{n}_K^H : X^K \to X^H$, defined by $n_K^H(x) = r_1 x * r_2 x * \cdots * r_{n} x$ for some choice of $H/K$-coset representatives $r_1,\dots,r_{n}$. Indeed, if $x \in X^K$ and $h \in H$, and we write $h \cdot r_i K = r_{\sigma i} K$, then
	\begin{align*}
	h \cdot (r_1  x * \cdots * r_{n} x ) &= hr_1 x * \cdots * hr_{n}x = r_{\sigma 1}x * \cdots * r_{\sigma n} x	= r_1x * \cdots * r_{n} x
	\end{align*}
by the strict $G$-equivariance and commutativity of $*$. Thus, the fixed-point presheaf of $X$ is a topological semi-Mackey functor.

While strict associativity and unitality are negotiable in homotopical algebra, strict commutativity is far too much to ask for. We say that a $G$-operad $\s{O}$ is \emph{$\Sigma$-free} if the $\Sigma_n$-action on $\s{O}(n)$ is free for every $n \geq 0$. Such operads parametrize no strict commutativity relations, and they typically have the most interesting algebras. 

\begin{ex} Suppose $V$ is a finite-dimensional real $G$-representation and write $D(V)$ for the unit disc centered at the origin in $V$. A \emph{little $V$-disc} in $D(V)$ is an affine, but not necessarily equivariant, map of the form $av + b : D(V) \to D(V)$. The $n$th level of the \emph{little $V$-discs operad} $\c{D}(V)$ is the space of all disjoint $n$-tuples of little $V$-discs in $D(V)$. The group $G$ acts on $\c{D}(V)(n)$ by conjugation, the group $\Sigma_n$ acts by permuting tuples, the map $\t{id} : D(V) \to D(V)$ is the operadic identity, and operadic composites are computed by slotting little $V$-discs into little $V$-discs. The operad $\c{D}(V)$ is $\Sigma$-free.
\end{ex}

The prototypical example of a $\c{D}(V)$-algebra $G$-space is the $V$-fold loop space $\Omega^V X = \t{Map}_*(S^V,X)$. Here $S^V$ is the one-point compactification of $V$, $X$ is a based $G$-space, and $\t{Map}_*(S^V,X)$ is the space of all continuous, based maps $S^V \to X$, equipped with the conjugation $G$-action. Conversely, every $\c{D}(V)$-algebra $G$-space group completes to a $V$-fold loop space, provided that $\bb{R}^2 \subset V$ \cite{GM}.

Experience has shown that transfer maps are useful and ubiquitous in genuine equivariant homotopy theory. However, they do not arise from the recipe above, because we very rarely have strictly commutative operations. For example, suppose $\rho$ is the regular representation of $G \neq \{e\}$. Then the sum on $\Omega^\rho X$ is only homotopy commutative. Nevertheless, for every pair of subgroups $K \subset H$, there is an additive transfer map $(\Omega^{\rho} X)^K \to (\Omega^\rho X)^H$. It arises by summing the $H$-conjugates of a $\rho$-loop $l \in (\Omega^{\rho}X)^K$ over a tubular neighborhood of a copy of $H/K \subset \t{res}^G_H \rho$. After ordering the orbit $H/K$, this neighborhood corresponds to an element $d \in \c{D}(\rho)(\abs{H:K})$, and this element is $H$-fixed, up to a twist given by the action of $H$ on $H/K$. We can formalize this kind of twisted equivariance, but first, a preliminary.

\begin{defn}Let $n \geq 0$ be a nonnegative integer. A \emph{graph subgroup} of $G \times \Sigma_n$ is a subgroup $\Gamma \subset G \times \Sigma_n$ that intersects $\Sigma_n = \{e\} \times \Sigma_n$ trivially.
\end{defn}

Crucially, if $\s{O}$ is a $\Sigma$-free $G$-operad and $f \in \s{O}(n)$, then $\t{Stab}(f) \subset G \times \Sigma_n$ is a graph subgroup. The terminology is motivated by the following standard observation.

\begin{lem}For any graph subgroup $\Gamma \subset G \times \Sigma_n$, there is a unique subgroup $H \subset G$ and group homomorphism $\sigma : H \to \Sigma_n$ such that $\Gamma = \{ (h,\sigma(h)) \, | \, h \in H\}$. Conversely, every subgroup of the form $\{(h,\sigma(h)) \, | \, h \in H\}$ is a graph subgroup.
\end{lem}

Now suppose that $\Gamma = \{ (h,\sigma(h)) \, | \, h \in H \} \subset G \times \Sigma_n$ is a graph subgroup, and that $f \in \s{O}(n)$ is a $\Gamma$-fixed operation. Then for any $\s{O}$-algebra $G$-space $X$, we obtain a $n$-ary product $f : X^{\times n} \to X$ such that
	\[
	hf(x_1,\dots,x_n) = f(hx_{\sigma(h)^{-1}1},\dots, h x_{\sigma(h)^{-1}n}) 
	\]
for every $h \in H$ and $(x_1,\dots,x_n) \in X^{\times n}$. Assume further that $\sigma : H \to \Sigma_n$ represents the $H$-action on $H/K = \{r_1 K < \dots < r_n K\}$, i.e. $h r_i K = r_{\sigma(h)i} K$ for every $h \in H$ and $1 \leq i \leq n$. Thus $f$ exhibits precisely the same equivariance as the operation $d \in \c{D}(\rho)$ considered above, and we obtain a norm map
	\[
	\t{n}_K^H(x) = f(r_1 x, \dots , r_n x) : X^K \to X^H .
	\]
Thus, if we are interested in constructing transfer maps in homotopy commutative settings, then a system of twisted equivariant maps, such as $f$ above, is a reasonable substitute for a strictly $G$-equivariant and commutative product $* : X^{\times 2} \to X$. Accordingly, we introduce the following terminology.

\begin{defn}\label{defn:extnorm} Suppose $X$ is a $G$-space, $H \subset G$ is a subgroup, and $T$ is a finite, ordered $H$-set whose permutation representation is $\sigma : H \to \Sigma_{\abs{T}}$. Write $\Gamma(T) = \{ (h,\sigma(h)) \, | \, h \in H\} \subset G \times \Sigma_n$ for the graph of $\sigma$. An \emph{external $T$-norm on $X$} is a $\Gamma(T)$-fixed point of $\b{End}(X)(\abs{T})$.
\end{defn}

Alternatively, if $X^{\times T}$ is the $T$-indexed power of $X$, i.e. the space $X^{\times \abs{T}}$ equipped with the $H$-action 
	\[
	h(x_1,\dots , x_{\abs{T}}) = (h x_{\sigma(h)^{-1}1} , \dots , h x_{\sigma(h)^{-1}\abs{T}}),
	\]
then an external $T$-norm on $X$ is an $H$-equivariant map $f : X^{\times T} \to X$.

With these notions in mind, we introduce $N_\infty$ operads. 

\begin{defn}\label{defn:Ninftyop} Let $\s{O}$ be a symmetric operad in the category $\b{Top}^G$ of $G$-spaces. We say that $\s{O}$ is a \emph{$N_\infty$ operad} if it satisfies the following three conditions:
	\begin{enumerate}
		\item{}for every integer $n \geq 0$, the $G \times \Sigma_n$-space $\s{O}(n)$ is $\Sigma_n$-free,
		\item{}for every graph subgroup $\Gamma \subset G \times \Sigma_n$, the subspace $\s{O}(n)^\Gamma$ is either empty or contractible, and
		\item{}the spaces $\s{O}(0)^G$ and $\s{O}(2)^G$ are nonempty.
	\end{enumerate}
We write $N_\infty\t{-}\b{Op}^G$ for the category of all $N_\infty$ $G$-operads.
\end{defn}

\begin{rem}This is equivalent to \cite[Definition 3.7]{BH1}. Note that (1) implies $\s{O}(n)^{\Xi} = \varnothing$ for all non-graph subgroups $\Xi \subset G \times \Sigma_n$, and that (3) implies $\s{O}(n)^G \neq \varnothing$ for all $n \geq 0$ because $\s{O}$ is a $G$-operad. Therefore $\s{O}(n)$ is a universal space for a family of subgroups of $G \times \Sigma_n$, which contains $H \times \{1\}$ for all subgroups $H \subset G$. In particular, $\s{O}(0)$ and $\s{O}(1)$ are $G$-contractible.
\end{rem}

Condition (2) ensures that $\s{O}$ parametrizes at most one external $T$-norm of each kind, up to coherent homotopy, and condition (3) ensures that $\s{O}$ parametrizes a homotopy coherent associative, commutative, and unital operation, for which all data is $G$-equivariant. More precisely, the $G$-fixed suboperad $\s{O}^G \subset \s{O}$ is an $E_\infty$ operad in the nonequivariant sense. Informally, we think of $N_\infty$ operads as representing objects for homotopy coherent incomplete semi-Mackey functors. A $N_\infty$ operad $\s{O}$ such that $\s{O}(n)^{\Gamma} \simeq *$ for every graph subgroup $\Gamma$ is often called an \emph{$E_\infty$ $G$-operad} (e.g. in \cite{LMS}, \cite{CostWan}, and \cite{GM}).

\begin{ex} Let $U$ be a $G$-universe, i.e., a countably infinite-dimensional real $G$-inner product space that contains each of its finite-dimensional subrepresentations infinitely often, and which also contains trivial summands.

The $n$th level of the \emph{linear isometries operad} $\c{L}(U)$ is the space of all linear, but not necessarily equivariant, isometries $U^{\oplus n} \to U$. The operad structure is inherited from $\b{End}(U)$. The operad $\c{L}(U)$ is $N_\infty$, and we think of it as representing the canonical multiplicative structure for $G$-spectra indexed over $U$.

The \emph{infinite little discs operad} $\c{D}(U)$ is the colimit $\t{colim}_{V \subset U} \c{D}(V)$ of the little $V$-discs operads $\c{D}(V)$, as $V$ ranges over all finite-dimensional subrepresentations of $U$. The operad $\c{D}(U)$ is $N_\infty$, and we think of it as representing the canonical additive structure for $G$-spectra indexed over $U$. However, there is a catch. The point-set level colimit that defines $\c{D}(U)$ is not compatible with suspension, and therefore $\c{D}(U)$ does not naturally act on infinite loop spaces structured by $U$. One can replace $\c{D}(U)$ with a levelwise homotopy-equivalent operad $\c{K}(U)$, called the \emph{infinite Steiner operad}, which does act on equivariant infinite loop spaces \cite{GM}.

Surprisingly, there are universes $U$ such that the operads $\c{D}(U)$ and $\c{L}(U)$ are inequivalent \cite[Theorem 4.22]{BH1}.
\end{ex}

Algebras over $N_\infty$ operads also appear in equivariant homotopical algebra for conceptual reasons. Hill and Hopkins \cite{HH} have proven that localizations of genuine commutative ring $G$-spectra need not have all multiplicative norms.  The underlying multiplication survives for formal reasons, which further justifies condition (3) in Definition \ref{defn:Ninftyop}, but that is all we are guaranteed. Subsequent work of Guti\'{e}rrez and White \cite{GutWhite} addresses when general left Bousfield localizations preserve and destroy $N_\infty$ algebra structures.

\subsection{The homotopy theory of $N_\infty$ operads} The purpose of a $N_\infty$ operad is to parametrize homotopy coherent algebraic structures. Accordingly, we introduce the following weak equivalences.

\begin{defn} An operad map $\vp : \s{O}_1 \to \s{O}_2$ between $N_\infty$ operads is a \emph{weak equivalence} if $\vp_n : \s{O}_1(n)^\Gamma \to \s{O}_2(n)^\Gamma$ is a weak homotopy equivalence of topological spaces for every $n \geq 0$ and graph subgroup $\Gamma \subset G \times \Sigma_n$.
\end{defn}

Note that a weak equivalence $\vp : \s{O}_1 \to \s{O}_2$ between $N_\infty$ operads is actually a levelwise weak $G \times \Sigma_n$-homotopy equivalence, because we have no $\Xi$-fixed points when $\Xi \subset G \times \Sigma_n$ is not a graph subgroup.

In contrast to the situation fo nonequivariant $E_\infty$ operads, not all $N_\infty$ $G$-operads are equivalent. However, May's product trick still works.

\begin{lem}\label{lem:prodtrick} Let $\s{O}_1$ and $\s{O}_2$ be $N_\infty$ operads. Suppose that for every $n \geq 0$ and graph subgroup $\Gamma \subset G \times \Sigma_n$, either $\s{O}_1(n)^\Gamma$ and $\s{O}_2(n)^\Gamma$ are both empty, or $\s{O}_1(n)^\Gamma$ and $\s{O}_2(n)^\Gamma$ are both nonempty. Then $\s{O}_1$ and $\s{O}_2$ are equivalent.
\end{lem}

\begin{proof} Both projections $\s{O}_1 \leftarrow \s{O}_1 \times \s{O}_2 \to \s{O}_2$ are weak equivalences.
\end{proof}

Thus, a $N_\infty$ operad $\s{O}$ is determined by the norms it parametrizes, or more formally, by the set of graph subgroups $\Gamma \subset G \times \Sigma_n$ such that $\s{O}(n)^\Gamma \neq \varnothing$. These collections cannot be arbitrary. If we fix $n \geq 0$, then the set of such $\Gamma$ is closed under subconjugacy. As we vary $n$, the operad structure on $\s{O}$ implies further closure conditions. It is convenient to phrase these conditions in terms of actions by subgroups of $G$.

\begin{defn}Suppose $\s{O}$ is a $N_\infty$ operad, $H \subset G$ is a subgroup, and $T$ is a finite $H$-set. Choose an order on $T$ and let $\Gamma(T) \subset G \times \Sigma_{\abs{T}}$ be the graph of the corresponding permutation representation. We say that \emph{$T$ is admissible for $\s{O}$} or that \emph{$\s{O}$ admits $T$} if $\s{O}(\abs{T})^{\Gamma(T)}$ is nonempty. We write $A(\s{O})$ for the class of all admissible sets of $\s{O}$.
\end{defn}

The admissibility of a $H$-set $T$ is independent of the choice of order on $T$ because different choices conjugate $\Gamma(T)$. Note that the class of admissible sets of a $N_\infty$ operad is graded over $\b{Sub}(G)$, the set of all subgroups of $G$.

\begin{defn}\label{defn:indsys}A \emph{class of finite $G$-subgroup actions} is a class $\c{X}$, equipped with a function $\c{X} \to \b{Sub}(G)$, such that the fiber over $H \subset G$ is a class of finite $H$-sets. We write $\c{X}(H)$ for the fiber over $H$. A \emph{$G$-indexing system} is a class of finite $G$-subgroup actions $\c{I}$ that satisfies the following seven conditions: 
	\begin{enumerate}
		\item{}(trivial sets) For any subgroup $H \subset G$, the class $\c{I}(H)$ contains all finite trivial $H$-actions.
		\item{}(isomorphism) For any subgroup $H \subset G$ and finite $H$-sets $S$ and $T$, if $S \in \c{I}(H)$ and $S \cong T$, then $T \in \c{I}(H)$.
		\item{}(restriction) For any subgroups $K \subset H \subset G$ and finite $H$-set $T$, if $T \in \c{I}(H)$, then $\t{res}^H_K T \in \c{I}(K)$.
		\item{}(conjugation) For any subgroup $H \subset G$, group element $a \in G$, and finite $H$-set $T$, if $T \in \c{I}(H)$, then $c_a T \in \c{I}(aHa^{-1})$.
		\item{}(subobjects) For any subgroup $H \subset G$ and finite $H$-sets $S$ and $T$, if $T \in \c{I}(H)$ and $S \subset T$, then $S \in \c{I}(H)$.
		\item{}(coproducts) For any subgroup $H \subset G$ and finite $H$-sets $S$ and $T$, if $S \in \c{I}(H)$ and $T \in \c{I}(H)$, then $S \sqcup T \in \c{I}(H)$.
		\item{}(self-induction) For any subgroups $K \subset H \subset G$ and finite $K$-set $T$, if $T \in \c{I}(K)$ and $H/K \in \c{I}(H)$, then $\t{ind}_K^H T \in \c{I}(H)$.
	\end{enumerate}
We call the elements of $\c{I}(H)$ the \emph{admissible $H$-sets of $\c{I}$}. Let $\b{Ind}(G)$ denote the class of all $G$-indexing systems.
\end{defn}

Condition (1) says the space $\s{O}(n)^G$ is nonempty for every $n \geq 0$. Conditions (2)--(4) say the set $\{ \Gamma \subset G \times \Sigma_n \, | \, \s{O}(n)^\Gamma \neq \varnothing \}$ is a family. Conditions (5)--(7) encode the operad structure on $\s{O}$. For every $k , j_1 , \dots , j_k \geq 0$, we have a $G$-equivariant composition map $\gamma : \s{O}(k) \times \s{O}(j_1) \times \dots \times \s{O}(j_k) \to \s{O}(j_1 + \cdots + j_k)$ that is also suitably $\Sigma$-equivariant. If the domain has a $\Gamma$-fixed point, then so does the codomain, and one can deduce conditions (5)--(7) by substituting specific norms into the domain.

\begin{rem}Indexing systems in the sense of Definition \ref{defn:indsys} are equivalent to indexing systems in the sense of \cite[Definition 3.22]{BH1}, because full subcategories are determined by their objects, and the axioms in Definition \ref{defn:indsys} imply closure under cartesian products. For suppose $S,T \in \c{I}(H)$ and choose orbit decompositions $S \cong \coprod_i H/K_i$ and $T \cong \coprod_j H/L_j$. Then $H/K_i, H/L_j \in \c{I}(H)$ for every $i$ and $j$ by (5), and $S \times T \cong \coprod_{i,j} (H/K_i \times H/L_j)$. By (2) and (6), it will be enough to show $H/K \times H/L \in \c{I}(H)$ whenever both $H/K \in \c{I}(H)$ and $H/L \in \c{I}(H)$, but this follows from the isomorphism $H/K \times H/L \cong \t{ind}_K^H \t{res}^H_K H/L$ and (2), (3), and (7).
\end{rem}

As suggested by the repeated use of ``admissible,'' we have the following result.

\begin{thm}[\protect{\cite[Theorem 4.17]{BH1}}] If $\s{O}$ is a $N_\infty$ $G$-operad, then the class $A(\s{O})$ of admissible sets of $\s{O}$ is a $G$-indexing system.
\end{thm}

This theorem is the key link between $N_\infty$ operads and indexing systems. Accordingly, we pause for a moment to analyze indexing systems.

If $\c{I}$ is a $G$-indexing system, then conditions (5) and (6) imply that $\c{I}(H)$ is the class of all finite coproducts of admissible $H$-orbits of $\c{I}$. Thus, $\c{I}$ is determined by the orbits it contains, and there are only finitely many $G$-indexing systems for a given group $G$.

Next, we declare $\c{I} \leq \c{J}$ if $\c{I}(H) \subset \c{J}(H)$ for every subgroup $H \subset G$. The componentwise intersection of a set of $G$-indexing systems is a $G$-indexing system, and therefore $(\c{I} \land \c{J})(H) = \c{I}(H) \cap \c{J}(H)$ is the meet of $\c{I}$ and $\c{J}$ in $\b{Ind}(G)$. The componentwise union $(\c{I} \cup \c{J})(H) = \c{I}(H) \cup \c{J}(H)$ is not always an indexing system, but it generates one.

\begin{defn}\label{defn:genind} For any class of finite $G$-subgroup actions $\c{X}$, we define $\la \c{X} \ra$ to be the intersection of all $G$-indexing systems that contain $\c{X}$.
\end{defn}

The join $\c{I} \lor \c{J}$ of $\c{I}$ and $\c{J}$ is the indexing system $\la \c{I} \cup \c{J} \ra$. It follows that $\b{Ind}(G)$ is a finite lattice. There is a maximum $G$-indexing system, whose $H$-component contains all finite $H$-sets, and there is a minimum $G$-indexing system, whose $H$-component contains only trivial finite $H$-sets. We denote the former $\ub{Set}$ and the latter $\ub{triv}$. We summarize.

\begin{prop}The class $\b{Ind}(G)$ of all $G$-indexing systems is a finite lattice under levelwise inclusion. The meet of two indexing systems is their levelwise intersection, the join of two indexing systems is the indexing system generated by their levelwise union, the minimum indexing system $\ub{triv}$ is class of all trivial actions, and the maximum indexing system $\ub{Set}$ is the class of all actions.
\end{prop}

We return to the classification of $N_\infty$ operads. Taking admissible sets sends an $N_\infty$ operad $\s{O}$ to an indexing system $A(\s{O})$, and converts a map $\vp : \s{O}_1 \to \s{O}_2$ between $N_\infty$ operads into an inclusion $A(\s{O}_1) \subset A(\s{O}_2)$. Moreover, if $\vp$ is a weak equivalence, then $A(\s{O}_1) = A(\s{O}_2)$. Thus we obtain a functor
	\[
	A : \t{Ho}(N_\infty\t{-}\b{Op}^G) \to \b{Ind}(G),
	\]
where $\t{Ho}(N_\infty \t{-}\b{Op}^G)$ is the category of $N_\infty$ operads with weak equivalences inverted. The classification theorem says this functor is an equivalence.

To show $A : \t{Ho}(N_\infty\t{-}\b{Op}^G) \to \b{Ind}(G)$ is full, note $A(\s{O}_1 \times \s{O}_2) = A(\s{O}_1) \land A(\s{O}_2)$ for any $N_\infty$ operads $\s{O}_1$ and $\s{O}_2$. Thus, if $A(\s{O}_1) \subset A(\s{O}_2)$, then $A(\s{O}_1 \times \s{O}_2) = A(\s{O}_1)$ and the left projection map in $\s{O}_1 \leftarrow \s{O}_1 \times \s{O}_2 \to \s{O}_2$ is an equivalence. Therefore this zig-zag determines a morphism $\s{O}_1 \to \s{O}_2$ in $\t{Ho}(N_\infty\t{-}\b{Op}^G)$, which maps to $A(\s{O}_1) \subset A(\s{O}_2)$ in $\b{Ind}(G)$.

Establishing faithfulness is more involved. Blumberg and Hill proved that every derived mapping space $\t{Map}(\s{O}_1,\s{O}_2)$ in the hammock localization $L^H(N_\infty\t{-}\b{Op}^G)$ is either empty or contractible \cite[Proposition 5.5]{BH1}. The strategy is to resolve $\s{O}_1$ by free operads, and then to use the free-forgetful adjunction and the emptiness or contractibility of $\s{O}_2$'s fixed point subspaces. Taking connected components of $L^H(N_\infty\t{-}\b{Op}^G)$ shows that every hom set in $\t{Ho}(N_\infty\t{-}\b{Op}^G)$ is either empty or a point, so the functor $A : \t{Ho}(N_\infty\t{-}\b{Op}^G) \to \b{Ind}(G)$ cannot help but be faithful. We give a new proof of this result in \S\ref{sec:modstr} (cf. Corollary \ref{cor:LHNop}).

Lastly, Blumberg and Hill made the following conjecture. 

\begin{conj}\label{conj:BH} The functor  $A : \t{Ho}(N_\infty\t{-}\b{Op}^G) \to \b{Ind}(G)$ is surjective.
\end{conj}

This has since been proven. We show that the functor $A$ is surjective in \S\ref{sec:realize} and \S\ref{sec:assocNop}, and both Bonventre-Pereira \cite{BonPer} and Guti\'{e}rrez-White \cite{GutWhite} have given independent proofs. Our approaches are rather different. As explained in \S\ref{sec:intro}, each has its own set of advantages, and each highlights distinct features of $N_\infty$ theory. However, there is a common theme in our solutions, which we describe in \S\ref{subsec:compare}.

We arrive at the following conclusion.

\begin{thm}[Classification of $N_\infty$ operads]\label{thm:classNinfty} Taking admissible sets determines a Dwyer-Kan equivalence $A : L^H(N_\infty \t{-} \b{Op}^G) \to \b{Ind}(G)$ of simplicial categories and an ordinary equivalence $A : \t{Ho}(N_\infty \t{-} \b{Op}^G) \to \b{Ind}(G)$ of $1$-categories.
\end{thm}

\begin{proof}Combine \cite[Theorem 3.24]{BH1} with Theorems \ref{thm:freerealize} or \ref{thm:assocNop} of this paper, or the results in \cite{BonPer} or \cite{GutWhite}.
\end{proof}

\begin{rem} Indexing systems are a natural device for studying $N_\infty$ operads, but there are other equivalent and useful formulations.

When contemplating incomplete Tambara functors, it is convenient to think in terms of polynomial bispans in the category $\b{Set}^G_{fin}$ of finite $G$-sets, whose multiplicative legs are restricted to a subcategory $\s{D} \subset \b{Set}^G_{fin}$. This subcategory $\s{D}$ should be wide, pullback stable, and finite coproduct complete to ensure that the corresponding category of bispans is sensible. Blumberg and Hill prove that such \emph{indexing categories} $\s{D}$ are in bijective correspondence with indexing systems \cite{BH2}.

One can also recast the definition of an indexing system purely in terms of orbits, and the result is what we call an \emph{transfer system}. More precisely, a transfer system is a partial order on $\b{Sub}(G)$ that refines inclusion, and which is closed under conjugation and restriction. Transfer systems are useful in combinatorially intensive situations, and we prove that transfer systems and indexing systems are equivalent in \cite{RubGeomInd}. This notion was also discovered in striking, independent work of Balchin, Barnes, and Roitzheim \cite{BBR}, in which they prove that the lattices $\b{Ind}(C_{p^n})$ are isomorphic to associahedra.
\end{rem}

\section{Discrete $N$ operads}\label{sec:Nop}

In this section, we explain how to reduce problems about $N_\infty$ operads to discrete combinatorics. The key point is that $N_\infty$ operads contain no higher homotopical information. We leverage this to give a quick construction of $N_\infty$ operads from operads in $\b{Set}^G$ that have the same isotropy properties. We call these combinatorial objects $N$ operads, and we show that $N$ operads are equivalent to $N_\infty$ operads for all homotopical purposes (Theorem \ref{thm:Gsetmod}). We conclude with a few examples of $N$ operads that elaborate on Guillou and May's constructions \cite{GM}.

\subsection{$N$ operads}\label{subsec:Nop} Consider the following discrete analogue to a $N_\infty$ operad.

\begin{defn}\label{defn:NOp} Let $\s{O}$ be a symmetric operad in the category $\b{Set}^G$ of $G$-sets with respect to the cartesian product. We say that $\s{O}$ is a \emph{$N$ operad} if it satisfies the following two conditions:
	\begin{enumerate}
		\item{}for every integer $n \geq 0$, the $G \times \Sigma_n$-set $\s{O}(n)$ is $\Sigma_n$-free,
		\item{}the sets $\s{O}(0)^G$ and $\s{O}(2)^G$ are nonempty.
	\end{enumerate}
We write $N\t{-}\b{Op}^G$ for the category of $N$ operads in $\b{Set}^G$.

For any subgroup $H \subset G$ and finite $H$-set $T$, we say that $T$ is \emph{admissible} for $\s{O}$ or that $\s{O}$ \emph{admits} $T$ if the set $\s{O}(\abs{T})^{\Gamma(T)}$ is nonempty. We write $A(\s{O})$ for the class of admissible sets of $\s{O}$.
\end{defn}

We construct $N_\infty$ operads from $N$ operads by attaching cells to kill all homotopy. This must be done somewhat carefully to ensure that the end result is still an operad. We borrow a trick from \cite{GMM}.

Let $\b{Cat}$ be the category of small categories. The functor
	\[
	\t{Ob} : \b{Cat} \to \b{Set},
	\]
which sends a small category $\s{C}$ to its set of objects, has a right adjoint 
	\[
	\b{Cat} \leftarrow \b{Set} : \til{(-)}.
	\]
For any $X \in \b{Set}$, the category $\til{X}$ has object set $X$, and a unique morphism $(x,y) : x \to y$ for every pair $x,y \in X$. Therefore $\til{\varnothing} = \varnothing$ and $\til{X} \simeq *$ if $X \neq \varnothing$.

\begin{defn} Let $E : \b{Set} \to \b{Top}$ be the composite of $\til{(-)} : \b{Set} \to \b{Cat}$ with the classifying space functor $B : \b{Cat} \to \b{sSet} \to \b{Top}$.
\end{defn}

The functor $E$ preserves all finite limits because $B$ does and $\til{(-)}$ is a right adjoint. It follows immediately that $E$ induces a functor
	\[
	E : \b{Op}(\b{Set}^G) \to \b{Op}(\b{Top}^G)
	\]
between categories of operads. The next observation explains our notation.

\begin{lem}\label{lem:Btilfp} Suppose $X$ is a $G$-set, and let $\c{F} = \{H \subset G \, | \, X^H \neq \varnothing \}$. Then $EX$ is a universal space for the family $\c{F}$.
\end{lem}

\begin{proof} For any subgroup $H \subset G$, the functor $(-)^H$ is a finite limit because $G$ is a finite group. Therefore $(EX)^H \cong E(X^H)$, and this is empty if $X^H = \varnothing$, and contractible if $X^H \neq \varnothing$.
\end{proof}

To go the other way, we ignore topology.

\begin{defn} Let $(-)^u : \b{Top} \to \b{Set}$ be the forgetful functor.
\end{defn}

The functor $(-)^u$ also preserves all (finite) limits, so it induces a functor
	\[
	\b{Op}(\b{Set}^G) \leftarrow \b{Op}(\b{Top}^G) : (-)^u .
	\]
The functors $E$ and $(-)^u$ form a tight link between $N_\infty$ operads and $N$ operads.

\begin{prop}\label{prop:ENop} Let $G$ be a finite group.
	\begin{enumerate}[label=(\roman*)]
		\item{}If $\s{O}$ is a $N$ operad in $\b{Set}^G$, then $E\s{O}$ is a $N_\infty$ operad in $\b{Top}^G$ with the same admissible sets.
		\item{}If $\s{O}$ is a $N_\infty$ operad in $\b{Top}^G$, then $\s{O}^u$ is a $N$ operad in $\b{Set}^G$ with the same admissible sets.

	\end{enumerate}
\end{prop}

\begin{proof}We begin with $(i)$. Suppose $\s{O}$ is a $N$ operad. We apply Lemma \ref{lem:Btilfp} repeatedly to verify the conditions in Definition \ref{defn:Ninftyop}. For (1), if $\Xi \subset \{e\} \times \Sigma_n$ is a nontrivial subgroup, then $\s{O}(n)^\Xi = \varnothing$ because $\s{O}(n)$ is $\Sigma_n$-free, and therefore $E{\s{O}}(n)^\Xi = \varnothing$ as well. Thus $E{\s{O}}(n)$ is a $\Sigma_n$-free space. Condition (3) follows from $\s{O}(0)^G ,\s{O}(2)^G \neq \varnothing$, and condition (2) is immediate from Lemma \ref{lem:Btilfp}. Therefore $E{\s{O}}$ is a $N_\infty$ operad, and for any graph subgroup $\Gamma \subset G \times \Sigma_n$, we know that $\s{O}(n)^\Gamma$ is nonempty if and only if $E{\s{O}}(n)^\Gamma$ is nonempty. Thus $\s{O}$ and $E{\s{O}}$ have the same admissible sets.

Claim $(ii)$ holds because the functor $(-)^u$ preserves $\Sigma$-freeness, emptiness, and nonemptiness. 
\end{proof}

Even though $N$ operads are discrete, we can equip the category of all $N$ operads with a perfectly good homotopy theory by creating weak equivalences along the functor $E : N\t{-}\b{Op}^G \to N_\infty\t{-}\b{Op}^G$.

\begin{defn} A morphism $f : \s{O}_1 \to \s{O}_2$ of $N$ operads is a \emph{weak equivalence} if $E{f} : E{\s{O}_1}(n)^{\Gamma} \to E{\s{O}_2}(n)^{\Gamma}$ is a weak homotopy equivalence of topological spaces for all $n \geq 0$ and graph subgroups $\Gamma \subset G \times \Sigma_n$.
\end{defn}

Since the fixed points $E{\s{O}}_i(n)^{\Gamma}$ are either empty or contractible for $i = 1,2$, saying $Ef : E\s{O}_1 \to E\s{O}_2$ is a weak equivalence is the same as saying $\s{O}_2(n)^{\Gamma} \neq \varnothing$ implies $\s{O}_1(n)^{\Gamma} \neq \varnothing$ for all $n$ and $\Gamma$. This is a purely combinatorial condition with little dependence on $f$; however, the existence of an operad map $f : \s{O}_1 \to \s{O}_2$ implies that if $\s{O}_2(n)$ has a $\Gamma$-fixed point, then some $\Gamma$-fixed point of $\s{O}_2(n)$ lifts along $f$ to a $\Gamma$-fixed point of $\s{O}_1(n)$, namely the image of an element $x \in \s{O}_1(n)^\Gamma$.

Thus, we have a homotopical category $N_\infty\t{-}\b{Op}^G$ of $N_\infty$ operads, and a homotopical category $N\t{-}\b{Op}^G$ of $N$ operads. The functor $E$ preserves weak equivalences by definition, and it is straightforward to show the functor $(-)^u$ also preserves weak equivalences. The interesting thing is that $E$ and $(-)^u$ induce an equivalence of homotopy theories.

\begin{thm} \label{thm:Gsetmod} The homotopical functors $E : N\t{-}\b{Op}^G \rightleftarrows N_\infty\t{-}\b{Op}^G : (-)^u$ preserve admissible sets and induce Dwyer-Kan equivalences between the hammock localizations of $N\t{-}\b{Op}^G$ and $N_\infty\t{-}\b{Op}^G$.
\end{thm}

\begin{proof} Proposition \ref{prop:ENop} handles the claim about admissibles. The remainder of the proof is another application of May's product trick \cite{MayGILS}. Let $\s{O}$ be a $N_\infty$ $G$-operad. Then $E(\s{O}^u)$ is a $N_\infty$ operad with the same admissible sets by Proposition \ref{prop:ENop}, and therefore both of the product projections
	\[
	\s{O} \leftarrow \s{O} \times E(\s{O}^u) \to E(\s{O}^u)
	\]
are weak equivalences. Therefore $E \circ (-)^u$ and the identity functor on $N_\infty\t{-}\b{Op}$ are connected through a zig-zag of natural weak equivalences. Similar reasoning shows that $(-)^u \circ E$ and the identity functor on $N\t{-}\b{Op}$ are connected through a zig-zag of natural weak equivalences. Therefore $E$ and $(-)^u$ induce Dwyer-Kan equivalences between the simplicial hammock localizations $L^H(N_\infty\t{-}\b{Op}^G)$ and $L^H(N\t{-}\b{Op}^G)$ (cf. \cite[Propositions 3.3 and 3.5]{DK}).
\end{proof}

Thus, there is no homotopical difference between topological $N_\infty$ operads and discrete $N$ operads.

\begin{rem} Blumberg and Hill prove that every hom space in $L^H(N_\infty\t{-}\b{Op}^G)$ is either empty or contractible (cf. \cite[Proposition 5.5]{BH1}), so the same is true for the hom spaces in $L^H(N\t{-}\b{Op}^G)$. We shall give a purely combinatorial argument for this fact in \S\ref{sec:modstr}, thus reproving Blumberg and Hill's result.
\end{rem}

\subsection{Examples of $N$ operads}\label{subsec:exNOp} We now describe a few examples of $N$ operads that build on the ideas in \cite{GM}. We begin with coinduced operads.

Suppose $X$ is a nonempty right $G$-set and $\s{O}$ is a $N$ operad in $\b{Set}$, i.e. $\s{O}$  is $\Sigma$-free and $\s{O}(0) , \s{O}(2) \neq \varnothing$. Then $\b{Set}(X,\s{O})$ is a $N$ $G$-operad. Moreover, if $T$ is a finite $H$-set, then
	\[
	\t{$\b{Set}(X,\s{O})$ admits $T$}	\quad\t{if and only if}\quad
	\begin{array}{c}
		\t{every $h \in H$ that fixes a point in $X$}	\\
		\t{acts as the identity on $T$.}
	\end{array}
	\]
Here are two extremal cases of this construction.

\begin{ex}\label{ex:coindAs} Suppose $\b{As}$ is the associativity operad. Its $n$-ary operations are $\b{As}(n) = \Sigma_n$, with $\Sigma_n$ acting on the right by group multiplication. Let $X = G$, with $G$ also acting on the right by group multiplication. Then $\s{O} = \b{Set}(G,\b{As})$ is a $N$ operad, and $A(\s{O}) = \ub{Set}$. Applying the right adjoint to $\t{Ob} : \b{Cat} \to \b{Set}$ yields an operad $\til{\s{O}}$, which is isomorphic to the operad $\s{P}_G$ considered in \cite{GM} and \cite{GMMO}.

On the other hand, if $X = *$, then the $N$ operad $\s{O} = \b{Set}(X,\b{As})$ is isomorphic to $\b{As}$ equipped with a trivial $G$-action. Therefore $A(\s{O}) = \ub{triv}$, and $\til{\s{O}}$ is the ordinary Barratt-Eccles operad $\s{P}$ equipped with a trivial $G$-action.
\end{ex}

Unfortunately, not every indexing system $\c{I}$ is of the form $A(\b{Set}(X,\b{As}))$.

\begin{cex}\label{cex:C4coind} Let $G = C_4$ and choose a generator $g \in G$. Let $H = \{e,g^2\}$ and let $\c{I}$ be the $C_4$-indexing system that contains all finite $H$-sets, but only trivial sets otherwise. Then $\c{I} \neq A(\b{Set}(X,\b{As}))$ for every nonempty right $G$-set $X$. Indeed, if $\c{I} \subset A(\b{Set}(X,\b{As}))$, then $\b{Set}(X,\b{As})$ admits $H/e$, and then since $g^2$ acts nontrivially on $H/e$, it follows $g^2$ cannot fix any element of $X$. Therefore $G$ acts freely on $X$ and $A(\b{Set}(X,\b{As})) = \ub{Set}$ properly contains $\c{I}$.
\end{cex}

Similarly, one might hope that every $G$-indexing system is realized by a suboperad of $\b{Set}(G,\b{As})$, because $A(\b{Set}(G,\b{As})) = \ub{Set}$ is the terminal indexing system. This is also false. Bonventre shows that if $G = C_2 \times C_2$, then the indexing system that contains all finite $C_2 \times 1$-sets, but only trivial actions otherwise, cannot be realized as a suboperad of $\s{P}_G$ \cite[Example B.2.1]{BonThesis}. The problem is that the elements of $\b{Set}(G,\b{As})$ are overcrowded.

We now consider a discrete variant of the isometries operad, following \cite[\S 7]{GM}.

\begin{defn} A \emph{discrete $G$-universe} is a countably infinite $G$-set $U$, which contains infinitely many copies of each orbit $G/H$ that embeds in $U$, and which also contains copies of $G/G$.
\end{defn}

The following is a generalization of Guillou and May's additive operad $\s{V}_G(U)$.

\begin{ex}\label{ex:dlinisom} Suppose $U$ is a discrete $G$-universe. The $n$th level of the operad $\c{L}_{\t{d}}(U)$ is the set of all injective, but not necessarily equivariant, functions $U^{\sqcup n} \hookrightarrow U$, where $U^{\sqcup n}$ is the $n$-fold coproduct of $U$. The group $G$ acts by conjugation, $\Sigma_n$ acts by permuting $U$ summands, the identity function $\t{id} : U \to U$ is the identity, and $\gamma(g;f_1,\dots,f_k) = g \circ (f_1 \sqcup \cdots \sqcup f_k)$ is composition. The $\Sigma_n$ actions are free, and $\c{L}_{\t{d}}(U)(n)^G \neq \varnothing$ because $U^{\sqcup n}$ $G$-embeds into $U$. Therefore $\c{L}_{\t{d}}(U)$ is a $N$ operad.
\end{ex}

The admissible sets of $\c{L}_{\t{d}}(U)$ are easy to calculate. Let $U$ be a discrete $G$-universe and for any $H \subset G$, define $\t{Stab}_H(U) = \{\t{Stab}_H(x) \, | \, x \in U\}$. Then for any subgroups $K \subset H \subset G$,
	\[
	\t{$\c{L}_{\t{d}}(U)$ admits $H/K$}	\quad\t{if and only if}\quad
	\t{Stab}_K(U) \subset \t{Stab}_H(U) .
	\]
Consequently, not every indexing system is realized by an operad $\c{L}_d(U)$.

\begin{cex} Let $G = C_4$ and keep notation as in Counterexample \ref{cex:C4coind}. Then the indexing system $\c{I}$ is not realized by the operad $\c{L}_{\t{d}}(U)$ for any discrete $G$-universe $U$. For suppose $\c{L}_{\t{d}}(U)$ admits $H/ e $. Then $\{e\} = \t{Stab}_{e}(U) \subset \t{Stab}_H(U)$, hence $U$ contains the free orbit $C_4/ e $, and hence $\t{Stab}_{e}(U) \subset \t{Stab}_G(U)$. Therefore $\c{L}_{\t{d}}(U)$ also admits $C_4/e$.
\end{cex}

The relationship between $\c{L}_{\t{d}}(U)$ and the topological linear isometries operad $\c{L}(\bb{R}[U])$ is delicate. The extremal cases are easy. If $U = [G/G]^{\sqcup \infty}$, then $A(\c{L}_{\t{d}}(U)) = \ub{triv} = A(\c{L}(\bb{R}[U]))$, and if $U$ contains all $G$-orbits, then $A(\c{L}_{\t{d}}(U)) = \ub{Set} = A(\c{L}(\bb{R}[U]))$. Things are less clear in between. If $U = [G/G \sqcup G/e]^{\sqcup \infty}$, then $\bb{R}[U]$ is a complete $G$-universe and $A(\c{L}(\bb{R}[U])) = \ub{Set}$. On the other hand,  $\c{L}_{\t{d}}(U)$ does not admit $G/H$ for any nontrivial, proper subgroup $H$ when $U = [G/G \sqcup G/e]^{\sqcup \infty}$.

There is also a multiplicative variant of $\c{L}_{\t{d}}(U)$, which generalizes Guillou and May's operad $\s{V}_G^{\times}(U)$, and which is trying to model a linear isometries operad based on the tensor powers of a universe. We shall not pursue it here.

\section{The realization problem}\label{sec:realize}

Despite the counterexamples in \S\ref{subsec:exNOp}, it is possible to realize every indexing system by a $N$ operad or a $N_\infty$ operad. In this section, we give the simplest general construction that we know (Theorem \ref{thm:freerealize}). The linchpin of our work is Theorem \ref{thm:admfree}, a calculation that is logically equivalent to Blumberg and Hill's indexing system conjecture (Proposition \ref{prop:keycalcBHconj}). Its proof is somewhat involved, so we defer it to \S\ref{subsec:pfkeycalc}. We shall give a more refined construction of associative $N_\infty$ operads in \S\ref{sec:assocNop}.

\subsection{The key calculation}\label{subsec:keycalc}

We analyze the universal examples of $N$ operads and indexing systems. By general considerations, there is a free-forgetful adjunction
	\[
	F : \b{Sym}(\b{Set}^G) \rightleftarrows \b{Op}(\b{Set}^G) : U
	\]
between the categories of symmetric sequences and operads in $G$-sets. The left adjoint sends a $G$-symmetric sequence $S$ to the free $G$-operad $F(S)$ that it generates. There is an analogous adjunction for indexing systems, and miraculously, taking admissible sets preserves the adjunction, provided the operads and symmetric sequences are suitably restricted.

\emph{This is a non-formal fact.} It hinges on a calculation of the fixed points of a free $G$-operad, which amounts to composing a left adjoint with a right adjoint.

We begin on the operadic side, by restricting attention to $N$ operads and to certain symmetric sequences that generate them.

\begin{defn}\label{defn:Nsym}Let $S$ be a symmetric sequence in the category $\b{Set}^G$ of $G$-sets. We say $S$ is a \emph{$N$ symmetric sequence} if:
	\begin{enumerate}
		\item{}for every integer $n \geq 0$, the $G \times \Sigma_n$-set $S(n)$ is $\Sigma_n$-free, and
		\item{}the sets $S(0)^G$ and $S(2)^G$ are nonempty.
	\end{enumerate}
We write $N\t{-}\b{Sym}^G$ for the category of all $N$ symmetric sequences in $\b{Set}^G$.

For any subgroup $H \subset G$ and finite $H$-set $T$, we say that $T$ is \emph{admissible} for $S$ if $S(\abs{T})^{\Gamma(T)}$ is nonempty.
\end{defn}

By neglect of structure, every $N$ operad $\s{O}$ is a $N$ symmetric sequence. Conversely, every $N$ symmetric sequence generates a $N$ operad.

\begin{prop}\label{prop:Nfreeforget} The free-forgetful adjunction $F : \b{Sym}(\b{Set}^G) \rightleftarrows \b{Op}(\b{Set}^G) : U$ restricts to an adjunction
	\[
	F : N\t{-}\b{Sym}^G \rightleftarrows N\t{-}\b{Op}^G : U.
	\]
between the full subcategories of $N$ symmetric sequences and $N$ operads.
\end{prop}

\begin{proof} It is enough to show that $F(S) \in N\t{-}\b{Op}$ for every $S \in N\t{-}\b{Sym}^G$. If $S \in N\t{-}\b{Sym}^G$, then there is an operad map $F(S) \to \b{Set}(G,\b{As})$. Therefore $F(S)$ is $\Sigma$-free, and $F(S)(n)^G \neq \varnothing$ for $n = 0,2$ because of the unit $\eta : S \to F(S)$.
\end{proof}

The admissible sets of a $N$ symmetric sequence do not form an indexing system, because the conditions on subobjects, coproducts, and self-induction reflect operadic composition. We do retain some of the axioms in Definition \ref{defn:indsys}, though.

\begin{defn}\label{defn:coeffsys} A class of $G$-subgroup actions $\c{X}$ is a \emph{$G$-coefficient system} if it satisfies conditions (2)--(4) of Definition \ref{defn:indsys}. Let $\b{Coef}(G)$ be the poset of all $G$-coefficient systems, ordered under inclusion.
\end{defn}

Coefficient systems in the sense above are equivalent to full, replete subcoefficient systems of $\ub{Set}$ in the sense of \cite{BH1}. Since the subgroups of $G \times \Sigma_n$ that have nonempty fixed points are closed under subconjugacy, the next result follows.

\begin{lem}\label{lem:admNsseq} If $Q$ is a $N$ symmetric sequence, then $A(Q)$ is a coefficient system.
\end{lem}

However, if $\s{O}$ is a $N$ operad, then we get all of the axioms. 

\begin{prop}\label{prop:admind} If $\s{O}$ is a $N$ operad, then $A(\s{O})$ is a indexing system.
\end{prop}

\begin{proof} We have $A(\s{O}) = A(E{\s{O}})$, and $E{\s{O}}$ is a $N_\infty$ operad. Alternatively, Blumberg and Hill's original arguments \cite[\S 4]{BH1} work just fine, once we replace all instances of ``contractible'' with ``nonempty.''
\end{proof}

We now turn to the analogue of $F : N\t{-}\b{Sym} \rightleftarrows N\t{-}\b{Op} : U$ for indexing systems. There is a free-forgetful adjunction
	\[
	\la \bullet \ra : \b{Coef}(G) \rightleftarrows \b{Ind}(G) : \iota
	\]
where $\iota$ is the inclusion, and $\la \bullet \ra$ sends a $G$-coefficient system to the indexing system that it generates (cf. Definition \ref{defn:genind}). Consider the squares below.
	\[
	\begin{tikzpicture}
		\node(Sym1) at (0,1.5) {$N\t{-}\b{Sym}$};
		\node(Op1) at (2,1.5) {$N\t{-}\b{Op}$};
		\node(Coef1) at (0,0) {$\b{Coef}(G)$};
		\node(Ind1) at (2,0) {$\b{Ind}(G)$};
		
		\node(Sym2) at (5,1.5) {$N\t{-}\b{Sym}$};
		\node(Op2) at (7,1.5) {$N\t{-}\b{Op}$};
		\node(Coef2) at (5,0) {$\b{Coef}(G)$};
		\node(Ind2) at (7,0) {$\b{Ind}(G)$};
		
		\path[->]
		(Sym1) edge [above] node {$F$} (Op1)
		(Coef1) edge [below] node {$\la \bullet \ra$} (Ind1)
		(Sym1) edge [left] node {$A$} (Coef1)
		(Op1) edge [right] node {$A$} (Ind1)
		
		(Op2) edge [above] node {$U$} (Sym2)
		(Ind2) edge [below] node {$\iota$} (Coef2)
		(Sym2) edge [left] node {$A$} (Coef2)
		(Op2) edge [right] node {$A$} (Ind2)
		;
	\end{tikzpicture}
	\]

For any $N$ operad $\s{O}$, the equality $A(U (\s{O})) = \iota ( A(\s{O}) )$ for right adjoints is immediate. The equality $A(F(S)) = \la A(S) \ra$ for left adjoints also holds, but this is the crux of the problem.

\begin{thm}\label{thm:admfree} If $S$ is a  $N$ symmetric sequence, then $A(F(S)) = \la A (S) \ra$.
\end{thm}

\begin{proof}[Sketch of Proof] The inclusion $\la A(S) \ra \subset A(F(S))$ follows from the equivariance of the unit $\eta : S \to F(S)$ and the fact that $A(F(S))$ is an indexing system. The other inclusion requires work. We unpack the general theory of generators and relations for operads in \S\ref{sec:quotfreeop}, and then we calculate the admissible sets of $F(S)$ in \S\ref{subsec:pfkeycalc}.
\end{proof}

\begin{rem}\label{rem:keycalc} Here is how to interpret the equality $A(F(S)) = \la A(S) \ra$. The indexing system $\la A(S) \ra$ is obtained from the symmetric sequence $S$ by taking the external norms of $S$ (cf. Definition \ref{defn:extnorm}), and then closing up under the indexing system axioms. On the other hand, the indexing system $A(F(S))$ is obtained by closing up $S$ under composition, and then computing the resulting external norms. In the former case, the closure conditions of an indexing system are dictated by Blumberg and Hill. In the latter case, the closure conditions on $A(F(S))$ are dictated by algebra. That $A(F(S))$ and $\la A(S) \ra$ are equal says that Blumberg and Hill's indexing system axioms perfectly capture the algebra of composition for external norms.
\end{rem}

\subsection{Free realizations of indexing systems}

Assuming Theorem \ref{thm:admfree}, we can construct operadic realizations of all indexing systems, thus verifying Blumberg and Hill's indexing system conjecture.

\begin{defn}\label{defn:FT} Let $\c{T} = (T_{\alpha})_{\alpha \in J}$ be an indexed set of finite $G$-subgroup actions, and let $S_{\c{T}}$ be the $N$ symmetric sequence
	\[
	S_{\c{T}} = \frac{G \times \Sigma_0}{G \times \{\t{id}_0\}} \,\, \sqcup \,\, \frac{G \times \Sigma_2}{G \times \{\t{id}_2\}} \,\, \sqcup \,\, \coprod_{\alpha \in J}  \frac{G \times \Sigma_{\abs{T_{\alpha}}} }{ \Gamma(T_{\alpha})} .
	\]
We define $\b{F}_{\c{T}} = F(S_{\c{T}})$ to be the free $N$-operad on $S_{\c{T}}$.
\end{defn}

\begin{thm}\label{thm:freerealize} The functors $A : N_\infty\t{-}\b{Op}^G \to \b{Ind}(G)$ and $A : N\t{-}\b{Op}^G \to \b{Ind}(G)$ have functorial sections. In particular, there is a section
	\[
	\b{F} : \b{Ind}(G) \to N\t{-}\b{Op}^G
	\]
given by the formula $\b{F}(\c{I}) = \b{F}_{\b{O}(\c{I})}$, where $\b{O}(\c{I})$ is the set of nontrivial orbits $H/K \in \c{I}$. The operad $\b{F}(\c{I})$ is a finitely generated free operad for every $\c{I} \in \b{Ind}(G)$.
\end{thm}

\begin{proof} By Theorem \ref{thm:admfree}, we have
	\[
	A(\b{F}(\c{I})) = \la A(S_{\b{O}(\c{I})}) \ra = \la \b{O}(\c{I}) \ra = \c{I}.
	\]
Therefore $\b{F}(\c{I})$ is a $N$ operad that realizes $\c{I}$, $E\b{F}(\c{I})$ is a $N_\infty$ operad that realizes $\c{I}$, and Conjecture \ref{conj:BH} is true. Moreover, if $\c{I} \subset \c{J}$, then $S_{\b{O}(\c{I})} \subset S_{\b{O}(\c{J})}$, and this inclusion induces a map $\b{F}(\c{I}) \to \b{F}(\c{J})$. Therefore $\b{F} : \b{Ind}(G) \to N\t{-}\b{Op}^G$ is a functorial section of $A : N\t{-}\b{Op}^G \to \b{Ind}(G)$ and $E \circ \b{F} : \b{Ind}(G) \to N_\infty\t{-}\b{Op}^G$ is a functorial section of $A : N_\infty\t{-}\b{Op}^G \to \b{Ind}(G)$.
\end{proof}

We use the set $\b{O}(\c{I})$ to generate the operad $\b{F}(\c{I})$ because it is efficient and reasonably canonical. Plenty of other choices are possible.

\begin{ex}\label{ex:choosegens} For each subgroup $H \subset G$, integer $n \geq 0$, and homomorphism $\sigma : H \to \Sigma_n$, we write $(n,\sigma)$ for the $H$-action on $\{1,\dots,n\}$ determined by $\sigma$. Given an arbitrary indexing system $\c{I}$, let $\b{N}(\c{I})$ be the set of all $(n,\sigma)$ contained in $\c{I}$. Then $\b{N}(\c{I})$ contains every admissible set of $\c{I}$ up to isomorphism, and we obtain a functorial section $\b{F}_{\b{N}(\c{I})} : \b{Ind}(G) \to N\t{-}\b{Op}^G$ of $A$.
\end{ex}

We conclude with a comment on the logical significance of Theorem \ref{thm:admfree}.

\begin{prop}\label{prop:keycalcBHconj} Theorem \ref{thm:admfree} is logically equivalent to Conjecture \ref{conj:BH}
\end{prop}

\begin{proof}The proof of Theorem \ref{thm:freerealize} shows that Theorem \ref{thm:admfree} implies Conjecture \ref{conj:BH}. Now suppose that Theorem \ref{thm:admfree} were false. Then there would be some $N$ symmetric sequence $S$ such that $\c{I} = \la A(S) \ra \subsetneq A(F(S))$. We claim that $\c{I}$ would be unrealizable. Suppose for contradiction that $\c{I} = A(\s{O})$ for some $N$ operad $\s{O}$. Then $A(S) \subset \c{I} = A(\s{O})$, and therefore there would be a map $S \to \s{O}$ of symmetric sequences. By adjunction, we would obtain an operad map $F(S) \to \s{O}$, and deduce
	\[
	A(\s{O}) = \c{I} = \la A(S) \ra \subsetneq A(F(S)) \subset A(\s{O}).\qedhere
	\]
\end{proof}

\section{Free and quotient $G$-operads}\label{sec:quotfreeop}

There are plenty of excellent treatments of operads in a symmetric monoidal category (e.g. \cite{Rezk} and \cite{Fresse}). There are also excellent discussions of combinatorial operads in $\b{Set}$ (cf. \cite{CCG} and \cite{Giraudo}). Unfortunately, we could not find an account of operads in $\b{Set}^G$ that met our needs. The proof of Theorem \ref{thm:admfree} hinges on delicate equivariant combinatorics, and we require an extremely precise description of free $N$ operads to carry it out. Thus, we shall spend this section building scaffolding.

The basic theory of combinatorial operads has many formal similarities to ordinary algebra, and we shall omit the most routine proofs. Unfortunately, this material is fairly dry. Therefore we begin by summarizing the relevant results in \S\ref{subsec:sumfree}, and then we flesh out the details in \S\S\ref{subsec:quotop}--\ref{subsec:opFS}. We recommend skimming or skipping the latter on a first reading.

\subsection{Summary}\label{subsec:sumfree} We give an explicit description of the free $N$ operad generated by a $N$ symmetric sequence. Suppose $S = (S(n))_{n \geq 0}$ is a $N$ symmetric sequence in $\b{Set}^G$. We think of the elements $f \in S_n$ as $n$-ary operations, and we will usually write them as functions $f(x_1,\dots,x_n)$ in $x_1,\dots,x_n$.

The free $N$ operad $F(S)$ is constructed from $S$ in two stages. First, we construct a $G$-operad $F_0(S)$, whose $n$-ary operations are formal composites of the operations in $S$, which contain each of the variables $x_1,\dots,x_n$ exactly once. For example, if $f \in S(3)$, $h \in S(2)$, $k \in S(1)$, and $\ell \in S(3)$, then the formal composites
	\[
	f(h(x_3,x_2),k(x_1),\ell(x_6,x_4,x_5)) \quad\t{and}\quad f(h(k(x_6),x_5),\ell(x_4,x_3,x_2),x_1)
	\]
are in $F_0(S)$.

Operadic composition $\gamma$ on $F_0(S)$ is defined by reindexing variables and then substituting functions into functions. For example,
	\[
	\gamma\Big( f(x_2,x_1,x_3) ; k(x_1) , h(x_2,x_1) , \ell(x_3,x_1,x_2) \Big) = f(h(x_3,x_2),k(x_1),\ell(x_6,x_4,x_5))
	\]
This requires a bit of explanation. The left hand side really is correct, because we want all arguments of $\gamma$ to be elements of $F_0(S)$. Now, the idea is to substitute $k(x_1)$ for $x_1$, $h(x_2,x_1)$ for $x_2$, and $\ell(x_3,x_1,x_2)$ for $x_3$ in $f(x_1,x_2,x_3)$, but this does not work because it produces something with three $x_1$'s. Therefore we replace $h(x_2,x_1)$ with $h(x_3,x_2)$ and $\ell(x_3,x_1,x_2)$ with $\ell(x_6,x_4,x_5)$ before substituting.

Now for the rest of the structure. The variable $x_1$ is the identity for $\gamma$. Right multiplication by a permutation $\sigma$ moves $x_i$ to $x_{\sigma(i)}$'s spot, e.g.
	\[
	\ell(x_2,x_1,x_3) \cdot (321) = \ell(x_3,x_2,x_1) ,
	\]
and the $G$-action on $F_0(S)$ is inherited from the $G$-action on $S$, e.g.
	\[
	g \cdot (f(k(x_1),h(x_3,x_4,x_2))) = gf(gk(x_1),gh(x_3,x_4,x_2))
	\]
for any $g \in G$. We think of the $G$-action as conjugation, which commutes with composition. There is a natural inclusion map $\eta_0 : S \to F_0(S)$, which sends $f \in S_n$ to $f(x_1,\dots,x_n) \in F_0(S)$. This is the unit of the free-forgetful adjunction 
	\[
	F_0 : (\b{Set}^G)^{\bb{N}} \rightleftarrows \b{Op}(\b{Set})^G : U
	\]
between non-symmetric sequences of $G$-sets and operads in $G$-sets.

Unfortunately, the map $\eta_0 : S \to F_0(S)$ is not $\Sigma$-equivariant, because we forgot the $\Sigma$-action on $S$ when we constructed $F_0(S)$. We fix this by passing to a quotient. The operad $F(S)$ is the operad $F_0(S)$, modulo the relations
	\begin{equation*}
	f\sigma(x_1,\dots,x_n) \sim f(x_{\sigma^{-1}1}, \dots , x_{\sigma^{-1}n}) \quad\quad(n \geq 0 \, , \, f \in S(n) \, , \, \sigma \in \Sigma_n ).
	\end{equation*}
We write $[t]$ for the congruence class of $t \in F_0(S)$. Combining the universal properties of the unit $\eta_0 : S \to F_0(S)$ and the quotient $\pi : F_0(S) \to F(S)$ shows that $F(S)$ is the free $N$ operad generated by $S$. The unit of the free-forgetful adjunction
	\[
	F : \b{Sym}(\b{Set}^G) \rightleftarrows \b{Op}(\b{Set}^G) : U
	\]
is the composite $\pi \circ \eta_0 : S \to F_0(S) \to F(S)$, which sends $f$ to $[f(x_1,\dots,x_n)]$.

We shall need a more precise description of $F(S)$ in order to compute its fixed points. Choose a set $S_\Sigma$ of $\Sigma$-orbit representatives of $S$. By stringing together the $\sim$ relations above, we may convert any formal composite of operations in $S$ into a formal composite of operations in $S_\Sigma$. For example, if $f = f' \cdot (12) \in S(2)$, $h = h' \cdot (123) \in S(3)$, and $k = k' \in S(1)$, where $f', \, h', \, k' \in S_\Sigma$, then
	\[
	f(k(x_1),h(x_3,x_4,x_2)) \sim f'(h'(x_2,x_3,x_4),k'(x_1)) .
	\]
It follows that $F_0(S_\Sigma) \subset F_0(S)$ is a set of representatives for $\sim$, which implies it inherits an operad structure from $F(S)$. All of the structure on $F_0(S_\Sigma)$ is the same as in $F_0(S)$, except for the $G$-action. This is because $F_0(S_\Sigma)$ is not closed under the $G$-action of $F_0(S)$. Thus, for any $g \in G$ and $t \in F_0(S_\Sigma)$, we define a new product $g * t$ by computing $g \cdot t$ in $F_0(S)$, and then applying $\sim$ relations to convert $g \cdot t$ into an element of $F_0(S_\Sigma)$. For example, if $f \in S_\Sigma(n)$ is $\Gamma(T)$-fixed for some $n$-element $G$-set $T$, then for any $(g,\sigma(g)) \in \Gamma(T)$, we have $g \cdot f = f \cdot \sigma(g)$ in $F_0(S)$. Therefore
	\[
	g * f(x_1,\dots,x_n) = f(x_{\sigma(g)^{-1}1},\dots,x_{\sigma(g)^{-1}n}) \quad\t{ in $F_0(S_\Sigma)$,}
	\]
which means that $f(x_1,\dots,x_n)$ is formally an external $T$-norm. The operad $F_0(S_\Sigma)$, equipped with $*$, is isomorphic to the free operad $F(S)$. This is the model of $F(S)$ that we will use in \S\ref{subsec:pfkeycalc}. 

We shall spend the remainder of this section making the sketch above precise. We treat formal composites in \S\ref{subsec:formcomp}, we construct $F_0(S)$ in \S\ref{subsec:opFS0} (cf. Construction \ref{const:F0} and Proposition \ref{prop:F0}), and we construct $F(S)$ in \S\ref{subsec:opFS} (cf. Construction \ref{const:freeop} and Theorem \ref{thm:freeGsigfree}). Quotients operads are discussed in \S\ref{subsec:quotop}, because they logically precede the construction of $F(S)$. 

\subsection{Quotient operads}\label{subsec:quotop} Suppose $\s{O}$ is an operad in $\b{Set}^G$. Since composition in $\s{O}$ is typically non-invertible, we cannot construct quotients of $\s{O}$ as sets of cosets, as one typically does in group, ring, and module theory. We shall use congruence relations instead. They should be thought of as substitutes for normal subgroups, ideals, and submodules.

\begin{defn}\label{defn:congrel} Suppose $\s{O}$ is an operad in $\b{Set}^G$. A \emph{congruence relation} $\sim$ on $\s{O}$ is a tuple $(\sim_n)_{n \geq 0}$ such that
	\begin{enumerate}
		\item{}for all integers $n \geq 0$, $\sim_n$ is an equivalence relation on $\s{O}(n)$,
		\item{}for all integers $n \geq 0$, elements $(g,\sigma) \in G \times \Sigma_n$, and operations $f,f' \in \s{O}(n)$, if $f \sim_n f'$, then $gf\sigma \sim_n gf'\sigma$, and
		\item{}for all integers $k,j_1,\dots,j_k \geq 0$ and operations $h,h' \in \s{O}(k)$, $f_1 , f_1' \in \s{O}(j_1)$, \dots, $f_k , f_k' \in \s{O}(j_k)$, if $h \sim_k h'$ and $f_i \sim_{j_i} f_i'$ for $i = 1, \dots , k$, then $\gamma(h;f_1,\dots,f_k) \sim_{j_1 + \cdots + j_k} \gamma(h' ; f_1', \dots , f_k')$.
	\end{enumerate}
\end{defn}
In other words, a congruence relation on a $G$-operad is a graded equivalence relation that is compatible with the operad structure. The axioms for a congruence relation ensure that all of the structure on $\s{O}$ descends to congruence classes.

\begin{defn}\label{defn:quotop} Suppose $\s{O}$ is an operad in $\b{Set}^G$ and $\sim$ is a congruence relation on $\s{O}$. The \emph{quotient operad} $\ol{\s{O}} = \s{O}/ \! \sim$ is defined as follows:
	\begin{enumerate}
		\item{}the set $\ol{\s{O}}(n)$ is the set $\s{O}(n)/\!\sim_n$ of all $\sim_n$-equivalence classes of $\s{O}(n)$,
		\item{}given a class $[f] \in \ol{\s{O}}(n)$ and $(g,\sigma) \in G \times \Sigma_n$, we define $g[f]\sigma := [gf\sigma]$,
		\item{}the identity of $\ol{\s{O}}$ is the class $[\t{id}]$ of the identity in $\s{O}$, and
		\item{}for any integers $k,j_1, \dots, j_k \geq 0$ and classes $[h] \in \ol{\s{O}}(k)$, $[f_1] \in \ol{\s{O}}(j_1)$, \dots, $[f_k] \in \ol{\s{O}}(j_k)$, we define $\gamma([h] ; [f_1] , \dots , [f_k]) := [ \gamma(h;f_1, \dots , f_k) ]$.
	\end{enumerate}
\end{defn}

Moreover, the quotient map $\pi : \s{O} \to \ol{\s{O}}$ has the usual universal property.

\begin{prop}\label{prop:UMPquot} Suppose $\s{O}$ is an operad in $\b{Set}^G$, $\sim$ is a congruence relation on $\s{O}$, and let $\pi : \s{O} \to \ol{\s{O}} = \s{O}/\!\sim$ be defined by $\pi(f) = [f]$. Then:
	\begin{enumerate}
		\item{}The map $\pi : \s{O} \to \ol{\s{O}}$ is an operad map, and $\pi(f) = \pi(f')$ whenever $f \sim f'$.
		\item{}If $\vp : \s{O} \to \s{O}'$ is an operad map such that $f \sim f'$ implies $\vp(f) = \vp(f')$, then there is a unique operad map $\ol{\vp} : \ol{\s{O}} \to \s{O}'$ such that $\vp = \ol{\vp} \circ \pi$.
	\end{enumerate}
\end{prop}

As one might hope, we can take quotients by kernels, but only after reinterpreting kernels as congruence relations.

\begin{defn}\label{defn:kernel} Suppose $\vp : \s{O} \to \s{O}'$ is a map of operads in $\b{Set}^G$. The \emph{kernel} of $\vp$ is the congruence relation $\sim_\vp = (\sim_{\vp,n})_{n \geq 0}$ on $\s{O}$, defined by
	\[
	f \sim_{\vp,n} f'	\quad\t{if and only if}\quad		\vp(f) = \vp(f')	\quad\quad( n \geq 0 \, , \, f,f' \in \s{O}(n)).
	\]
\end{defn}

Since $\sim$ is the kernel of $\pi : \s{O} \to \s{O}/\!\sim$, it follows that congruence relations on $\s{O}$ are the same thing as kernels of operad maps out of $\s{O}$.

Congruence relations are typically quite large, and in practice, we shall specify them using a small set of generators.

\begin{defn}Suppose $\s{O}$ is an operad in $\b{Set}^G$ and $R$ is a graded binary relation on $\s{O}$, i.e. $R = (R_n)_{n \geq 0}$ where $R_n$ is a binary relation on $\s{O}(n)$. Then the congruence relation generated by $R$ is
	\[
	\la R \ra_n = \Big\{ (f,f') \in \s{O}(n)^{\times 2} \, \Big| \, f \sim f' \t{ for all congruence relations $\sim \, \supset R$} \Big\},
	\]
i.e. $\la R \ra$ is the levelwise intersection of all congruence relations that contain $R$.
\end{defn}

The relation $\la R \ra$ is the smallest congruence relation that contains $R$. We introduce the relation $R$ into an operad $\s{O}$ by first enlarging $R$ to $\la R \ra$, and then taking the quotient $\s{O}/\la R \ra$. This quotient also has the expected universal property.

\begin{cor}\label{cor:UMPquot} Suppose $\s{O}$ is an operad in $\b{Set}^G$, $R$ is a graded binary relation on $\s{O}$, and let $\pi : \s{O} \to \s{O}/ \la R \ra$ be the quotient map. Then:
	\begin{enumerate}
		\item{}The map $\pi : \s{O} \to \s{O}/\la R \ra$ is an operad map, and $\pi(f) = \pi(f')$ if $fRf'$.
		\item{}If $\vp : \s{O} \to \s{O}'$ is an operad map such that $f R f'$ implies $\vp(f) = \vp(f')$, then there is a unique operad map $\ol{\vp} : \s{O}/\la R \ra \to \s{O}'$ such that $\vp = \ol{\vp} \circ \pi$.
	\end{enumerate}
\end{cor}

\begin{proof} Part (1) follows immediately from part (1) of Proposition \ref{prop:UMPquot}. For part (2), suppose $\vp : \s{O} \to \s{O}'$ is an operad map such that $fRf'$ implies $\vp(f) = \vp(f')$. Then $R$ refines $\t{ker}(\vp)$, and therefore $\la R \ra$ must, too. Thus, if $f \la R \ra f'$, then $\vp(f) = \vp(f')$, and the existence and uniqueness of $\ol{\vp} : \s{O}/\la R \ra \to \s{O}'$ follows from part (2) of Proposition \ref{prop:UMPquot}.
\end{proof}

It can be difficult to determine if two operations are identified by the congruence relation generated by $R$, but the following description of $\la R \ra$ can help.

\begin{defn}Suppose $R$ is a graded binary relation on an operad $\s{O}$ in $\b{Set}^G$. Given any integer $n \geq 0$ and operations $f_1, f_2 \in \s{O}(n)$, declare $f_1 \to f_2$ if 
	\[
	f_b = g \cdot \big( r \circ_k \gamma(s_b ; t_1, \dots , t_j) \big) \cdot \sigma	\quad (b=1,2)
	\]
for some
	\begin{enumerate}
		\item{}$(g, \sigma) \in G \times \Sigma_{n}$, and
		\item{} $r \in \s{O}(m)$, $s_1, s_2 \in \s{O}(j)$, and $t_a \in \s{O}(i_a)$ for $a = 1 , \dots , j$, 
	\end{enumerate}
such that $s_1 R_j s_2$ and the integers $n,m,k,j,i_1,\dots,i_j \geq 0$ satisfy  $1 \leq k \leq m$ and $n = m + i_1 + \cdots + i_j - 1$. Here $\circ_k$ denotes the $k$th partial composition product.
\end{defn}

The relation $\to$ is obtained by closing $R$ under the $G \times \Sigma$-action and certain composites. It is not usually an equivalence relation, so we generate one.

\begin{prop}\label{prop:gencong} Suppose $\s{O}$ is an operad in $\b{Set}^G$ and $R$ is a graded binary relation on $\s{O}$. Then $\la R \ra$ is the equivalence relation generated by $\to$ levelwise.
\end{prop}

\begin{proof} Let $\sim$ be the equivalence relation generated by $\to$. It is straightforward to check that $\sim$ refines every congruence relation that contains $R$. Thus, we only need to check that $\sim$ is a congruence relation. Consider $\to$. By construction, it satisfes (2) of Definition \ref{defn:congrel}, and it also satisfies a version of (3) where we only replace one of the operations $h$ or $f_i$. It follows that $\sim$ satisfies (2) and (3), and it is a graded equivalence relation by construction.
\end{proof}

\subsection{Formal composites}\label{subsec:formcomp}
 
We now turn to the constructions of $F_0(S)$ and $F(S)$, starting with a precise description of ``formal composites.'' We begin with some standard notions in formal logic.

\begin{defn}\label{defn:opterms} Suppose $S = (S(n))_{n \geq 0}$ is a sequence of $G$-sets. Regard the following formal symbols
	\[
	\begin{array}{rl}
		x_i	&	i = 1 , 2 , 3 , \dots \\
		f	&	f \in \coprod_{n \geq 0} S(n)	\\
		( \quad ) \quad ,	&	\t{(punctuation)}
	\end{array}
	\]
as the \emph{letters} in an alphabet $\Sigma(S)$. The elements of the free symmetric $G$-operad $F_0(S)$ will be suitable finite sequences of these letters.

A \emph{word} $w$ is a finite, ordered sequence $l_1 l_2 \cdots l_n$ of letters $l_i \in \Sigma(S)$. We write $\ve$ for the empty word. A \emph{subword} of $w$ is a word that is either $\ve$ or of the form $l_j l_{j+1} \cdots l_{k-1} l_k$ for some $1 \leq j \leq k \leq n$. The \emph{length} $\lambda(w)$ of the word $w = l_1 l_2 \cdots l_n$ is $n$, and $\lambda(\ve) = 0$.

A \emph{term} is any word constructed through the following recursion:
	\begin{enumerate}
		\item{}every variable $x_i$ is a term, and
		\item{}if $t_1 , \dots , t_n$ are terms and $f \in S(n)$, then $f(t_1,\dots,t_n)$ is also a term.
	\end{enumerate}
A \emph{subterm} of a word $w$ is a subword of $w$ that is also a term. The \emph{complexity} of a term $t$ is the length of the longest chain of nested pairs of left and right parentheses in $t$. For example,
	\[
	t = f(h(k(x_6),x_5),\ell(x_4,x_3,x_2),x_1)
	\]
has complexity $3$. Thus, if $t = f(t_1,\dots,t_n)$, then the complexity of each $t_i$ is strictly less than the complexity of $t$.

The \emph{arity} of a term $t$ is the number of distinct variable symbols $x_i$ that appear in $t$. We say that a $n$-ary term $t$ is \emph{operadic} if each of the variables $x_1 , \dots , x_n$ occur in $t$ exactly once.
\end{defn}

\begin{nota}\label{nota:tbar} Suppose $t$ is a term. We write $\ol{t}$ for the operadic term obtained from $t$ by reindexing the variables in $t$ as $x_1 , x_2 , \dots$ from left to right.
\end{nota}

\begin{ex} The unary term $s = p(x_1,x_1,x_1)$ is not operadic, and neither is the ternary term $t = p(x_2,x_3,x_4)$. However, we have
	\[
	\ol{s} = \ol{t} = p(x_1,x_2,x_3),
	\]
and both $\ol{s}$ and $\ol{t}$ are operadic.
\end{ex}

As we explain below, every term can be parsed into subterms, depending on the configuration of parentheses within it. Such decompositions can be interpreted as trees, but even though the corresponding pictures are intuitive, we find them unwieldy in calculations. Thus, we use the logical formalism instead.

\begin{defn} An \emph{initial segment} of a word $w = l_1 l_2 \cdots l_n$ is a word of the form $s = l_1 l_2 \cdots l_k$ for some $0 \leq k \leq n$. We understand $s = \ve$ if $k=0$, and we say $s$ is a \emph{strict} initial segment if $k < n$. Dually, a \emph{terminal segment} of $w$ is a word that is either $\ve$ or of the form $s = l_k l_{k+1} \cdots l_n$, and we say $s$ is \emph{strict} if $1 < k$.
\end{defn}

The key to parsing a term into subterms is the following parenthesis count. The proof of the following is a straightforward induction on complexity.

\begin{lem}\label{lem:termform}For any word $w$, write $L(w)$ for the number of left parentheses in $w$ and $R(w)$ for the number of right parentheses. Suppose that $t$ is a term. Then:
	\begin{enumerate}
		\item{}$L(t) = R(t)$.
		\item{}If $s$ is an initial segment of $t$, then $L(s) \geq R(s)$, and the inequality is strict if $2 \leq \lambda(s) < \lambda(t)$. In the latter case, $s$ is not a term.
		\item{}If $s$ is a terminal segment of $t$, then $L(s) \leq R(s)$, and the inequality is strict if $0 < \lambda(s) \leq \lambda(t) - 2$. In the latter case, $s$ is not a term.
	\end{enumerate}
\end{lem}

\begin{prop}\label{prop:parseterm} Suppose $m , n \geq 0$ are integers, $f \in S(m)$, $g \in S(n)$, and that $s_1 , \dots , s_m , t_1 , \dots , t_n$ are terms. If $f(s_1 , \dots , s_m) = g(t_1 , \dots , t_n)$ as words in $\Sigma(S)$, then $m = n$, $f = g$, and $s_i = t_i$ for $i = 1 , \dots , m$.
\end{prop}

\begin{proof}Suppose $f(s_1 , \dots , s_m) = g(t_1 , \dots , t_n)$. Then $f = g$ because they are the first letters. To show $s_1 = t_1$, it is enough to check that $s_1$ and $t_1$ have the same length. If $\lambda(s_1) < \lambda(t_1)$, then $s_1$ is a strict initial segment of $t_1$. Either $s_1$ is a variable and $t_1$ is not, or $\lambda(s_1) \geq 2$. The former case is clearly impossible, and the latter is ruled out by (2) of Lemma \ref{lem:termform}. Continue inductively. 
\end{proof}

Thus, it makes sense to speak of \emph{the} subterms of a given term.

\subsection{The operad $F_0(S)$}\label{subsec:opFS0}

We now construct the free $G$-operad $F_0(S)$ on a non-symmetric sequence of $G$-sets.

\begin{const}\label{const:F0} Let $S \in (\b{Set}^G)^{\bb{N}}$ be a sequence of $G$-sets, and define a symmetric operad $F_0(S)$ in $\b{Set}^G$ as follows.
	\begin{enumerate}
		\item{}Let $F_0(S)(n)$ be the set of all $n$-ary operadic terms in the alphabet $\Sigma(S)$. 
		\item{}Given $t \in F_0(S)(n)$ and $\sigma \in \Sigma_n$, let $t \cdot \sigma$ be the $n$-ary operadic term obtained from $t$ by replacing $x_i$ with $x_{\sigma^{-1}i}$ for each $i = 1 , \dots , n$. This makes $F_0(S)(n)$ into a right $\Sigma_n$-set.
		\item{}Given $g \in G$ we define a left $G$-action on all terms in $\Sigma(S)$ by the recursion:
			\begin{enumerate}
				\item{}$g \cdot x_n = x_n$ for $n = 1, 2 , 3, \dots$, and
				\item{}$g  \cdot f(t_1,\dots,t_n) = f' ( g \cdot t_1 , \dots , g \cdot t_n)$, where $f' = gf \in S$.
			\end{enumerate}
		This action multiplies every letter $f \in S$ in a term by $g$, and does nothing to the variables and punctuation. Therefore it restricts to a $G$-action on each set $F_0(S)(n)$, which commutes with the $\Sigma_n$-action.
		\item{}The identity element is $x_1$. It is $G$-fixed by definition.
		\item{}Given a $k$-ary operadic term $t$ and $j_i$-ary operadic terms $s_i$ for $i = 1 , \dots , k$, the composite $\gamma(t;s_1, \dots, s_k)$ is defined by:
			\begin{enumerate}
				\item{}adding $j_1 + \cdots + j_{i-1}$ to the subscript of every variable appearing $s_i$ -- call this new term $s_i'$ -- and then
				\item{}substituting the terms $s_1' , \dots , s_k'$ in for the variables $x_1 , \dots , x_k$ in $t$.
			\end{enumerate}
		These substitutions commute with the substitutions that define the $G$-action, and therefore $\gamma$ is $G$-equivariant.
	\end{enumerate}
There is a $G$-equivariant unit map $\eta_0 : S \to F_0(S)$ that sends the letter $f \in S(n)$ to the $n$-ary operadic term $f(x_1,\dots,x_n) \in F_0(S)(n)$. If $u \in S(0)$ we set $\eta_0(u) = u()$.
\end{const}

An important technical point is that every operadic term in $F_0(S)$ may be expressed canonically as a composite. Recall Notation \ref{nota:tbar}.

\begin{nota} Suppose $t \in F_0(S)$ and $t = f(t_1,\dots,t_n)$ for some $f \in S(n)$ and $j_i$-ary terms $t_i$. Then there is a unique $\sigma \in \Sigma_{j_1 + \cdots + j_n}$ such that
	\[
	f(t_1,\dots,t_n) = \gamma(\eta_0(f) ; \ol{t_1} , \dots , \ol{t_n}) \cdot \sigma .
	\]
We call the right hand side the \emph{standard decomposition} of $t$. If $u \in S(0)$, we understand the standard decomposition of $u()$ to be $\gamma(\eta_0(u); ) \cdot 1$.
\end{nota}

\begin{ex}The standard decomposition of $q(q(x_1,x_3),x_2)$ is
	\[
	q(q(x_1,x_3),x_2) = \gamma(q(x_1,x_2);q(x_1,x_2),x_1) \cdot (23) .
	\]
\end{ex}

With this decomposition in tow, we can establish the freeness of $F_0(S)$.

\begin{prop}\label{prop:F0} The map $\eta_0 : S \to F_0(S)$ in Construction \ref{const:F0} is the unit of the free-forgetful adjunction $F_0 : (\b{Set}^G)^{\bb{N}} \leftrightarrows \b{Op}(\b{Set}^G) : U$.
\end{prop}

\begin{proof} One checks that $\eta_0$ has the necessary universal property.

Suppose $\s{O}$ is an operad in $\b{Set}^G$ and $\vp : S \to \s{O}$ is a map of non-symmetric sequences. Then there is at most one operad map $\Phi : F_0(S) \to \s{O}$ that extends $\vp$ along $\eta_0 : S \to F_0(S)$. Indeed, let $t = \gamma(\eta_0(f) ; \ol{t_1} , \dots , \ol{t_n}) \cdot \sigma$ be the standard decomposition of $t$. Since $\Phi$ is an operad map, we must have $\Phi(t) = \gamma(\vp(f);\Phi(\ol{t_1}) , \dots, \Phi(\ol{t_n})) \cdot \sigma$ and $\Phi(x_1) = \t{id}$, which determines $\Phi$ recursively.

Now define $\Phi : F_0(S) \to \s{O}$ by the recursion above. Straightforward checks show that $\Phi$ is an operadic extension of $\vp$. For example, the standard decomposition of $f(x_1,\dots,x_n)$ is $\gamma(\eta_0(f);x_1,\dots,x_1) \cdot 1$, and hence
	\[
	\Phi(\eta_0(f)) = \Phi(f(x_1,\dots,x_n)) = \gamma(\vp(f);\Phi(x_1),\dots,\Phi(x_1)) \cdot 1 = \vp(f).
	\]
Therefore $\Phi$ extends $\vp$ along $\eta_0 : S \to F_0(S)$.

The map $\Phi$ preserves the identity by definition. 

To see that $\Phi$ is $\Sigma$-equivariant, note first that $\Phi(t \cdot \tau) = \Phi(t) \cdot \tau$ is automatic if $t = x_1$. If $t = f(t_1,\dots,t_n)$, then $t \cdot \tau = f(t_1',\dots,t_n')$ for some terms $t_i'$ such that $\ol{t_i} = \ol{t_i'}$. Thus, if $t = \gamma(\eta_0(f);\ol{t_1},\dots,\ol{t_n}) \cdot \sigma$, then the standard decomposition of $t \cdot \tau$ is $\gamma(\eta_0(f);\ol{t_1},\dots,\ol{t_n}) \cdot \sigma\tau$, and thus $\Phi(t \cdot \tau) = \gamma(\vp(f);\Phi(\ol{t_1}),\dots,\Phi(\ol{t_n})) \cdot \sigma\tau = \Phi(t) \cdot \tau$.

The rest of the proof is similar. One can induct on complexity to show $\Phi$ is $G$-equivariant, and $\Phi(\gamma(t;s_1,\dots,s_k)) = \gamma(\Phi(t);\Phi(s_1),\dots,\Phi(s_k))$ also follows by induction on the complexity of $t$.
\end{proof}

\subsection{The operad $F(S)$}\label{subsec:opFS} Finally, we construct the free operad $F(S)$ on a symmetric sequence of $G$-sets $S$.

\begin{const}\label{const:freeop} Suppose $S \in (\b{Set}^G)^\Sigma$ is a symmetric sequence of $G$-sets. Define the $G$-operad $F(S)$ by
	\[
	F(S) = \frac{F_0(S)}{\Big\la f\sigma(x_1,\dots,x_n) \sim f(x_{\sigma^{-1}1} , \dots , x_{\sigma^{-1}n}) \, \Big| \, n \geq 0 , \, f \in S(n) , \, \sigma \in \Sigma_n \Big\ra}
	\]
and let $\eta = \pi \circ \eta_0 : S \to F(S)$ be the composite of $\eta_0 : S \to F_0(S)$ and the projection $\pi : F_0(S) \to F(S)$, i.e. $\eta(f) = [f(x_1,\dots,x_n)]$.
\end{const}

\begin{prop}\label{prop:freequot} The map $\eta = \pi \circ \eta_0: S \to F(S)$ is the unit of the free-forgetful adjunction $F : (\b{Set}^G)^{\Sigma} \rightleftarrows \b{Op}(\b{Set}^G) : U$.
\end{prop}

\begin{proof}The relations that define the quotient $\pi : F_0(S) \to F(S)$ ensure that $\eta  = \pi \circ \eta_0$ is $G \times \Sigma$-equivariant, and universal property of $\eta$ follows from those of $\pi$ and of $\eta_0$ (cf. Proposition \ref{prop:F0} and Corollary \ref{cor:UMPquot}).
\end{proof}

We now give a more combinatorial description of $F(S)$ when $S$ is a $\Sigma$-free symmetric sequence of $G$-sets.

\begin{thm}\label{thm:freeGsigfree} Suppose that $S \in (\b{Set}^G)^{\Sigma}$ is a $\Sigma$-free symmetric sequence of $G$-sets and that $S_\Sigma(n) \subset S(n)$ is a set of $\Sigma_n$-orbit representatives for every integer $n \geq 0$. Then the free operad $F(S)$ on $S$ is isomorphic to the operad $F_0(S_\Sigma)$ in $\b{Set}$, equipped with the following recursively defined $G$-action. For any $g \in G$, declare:
	\begin{enumerate}
		\item{}$g * x_n = x_n$ for every $n>0$, and
		\item{}$g * f(t_1,\dots,t_n) = f'(g * t_{\sigma^{-1}1}, \dots , g * t_{\sigma^{-1}n})$, where $gf = f'\sigma$ for $f' \in S_\Sigma(n)$ and $\sigma \in \Sigma_n$, and the terms $t_1,\dots,t_n$ are not necessarily operadic.
	\end{enumerate}
If $f \in S(n)$, then the unit $\eta : S \to F_0(S_\Sigma)$ is defined by $\eta(f) = \eta_0(f')\sigma = f'(x_{\sigma^{-1}1},\dots,x_{\sigma^{-1}n})$, where $f = f'\sigma$ for $f' \in S_\Sigma(n)$ and $\sigma \in \Sigma_n$.
\end{thm}

\begin{proof} By Construction \ref{const:freeop}, the operad $F(S)$ is a quotient $F_0(S)/\!\!\sim$. By Proposition \ref{prop:gencong}, two $n$-ary terms $t,t'$ in $F_0(S)$ are identified by $\sim$ if and only if there is $m \geq 0$ and a sequence $t_0, \dots, t_m$ of $n$-ary terms of $F_0(S)$ such that 
	\begin{enumerate}
		\item{}$t = t_0$ and $t' = t_m$, and 
		\item{}for each $0 \leq i < m$, either the term $t_{i+1}$ is obtained by replacing a subterm of $t_i$ of the form $s = f\sigma(t_1,\dots,t_k)$ with the subterm $s' = f(t_{\sigma^{-1}1} , \dots , t_{\sigma^{-1}k})$, or vice versa.
	\end{enumerate}
We now give a simpler description of $\sim$. Declare $t \to t'$ if
	\begin{enumerate}[label=(\roman*)] 
		\item{}$t'$ is obtained by replacing a subterm of $t$ of the form $s = f\sigma(t_1,\dots,t_k)$ with the subterm $s' =  f(t_{\sigma^{-1}1} , \dots , t_{\sigma^{-1}k})$, and 
		\item{}$f \in S_\Sigma(k)$ and $\sigma \neq 1$.
	\end{enumerate}
Let $\stackrel{*}{\to}$ denote the reflexive and transitive closure of the relation $\to$. Then the diamond lemma implies that for every $t \in F_0(S)(n)$, there is a unique $r(t) \in F_0(S_\Sigma)(n)$ such that $t \stackrel{*}{\to} r(t)$.

If $t,t' \in F_0(S)(n)$ are such that $r(t) = r(t')$, then $t \sim t'$ becuase there is a a chain of forwards and backwards $\to$ relations between them. Conversely, if $t \sim t'$, then $r(t) = r(t')$. Indeed, it is enough to consider the case where $t = \alpha f\sigma(t_1,\dots,t_k)\beta$ and $t' = \alpha f(t_{\sigma^{-1}1},\dots,t_{\sigma^{-1}k})\beta$ for some words $\alpha$ and $\beta$, and $\sigma \neq 1$. If $f \in S_\Sigma$, then $t \to t'$, and hence $r(t) = r(t')$. If not, then $f = f_\Sigma \tau$ for $f_\Sigma \in S_\Sigma$ and $\tau \neq 1$. Writing $t'' = \alpha f_\Sigma(t_{(\tau\sigma)^{-1}1},\dots,t_{(\tau\sigma)^{-1}k})\beta$, we have $t \to t''$ and $t' \to t''$, so $r(t) = r(t'') = r(t')$.

It follows that $F_0(S_\Sigma)$ is a set of representatives for $\sim$ on $F_0(S)$, with $r(t)$ representing $t$. The quotient $\pi : F_0(S) \to F(S)$ induces a bijection $\pi : F_0(S_\Sigma) \to F(S)$ whose inverse is $\pi^{-1}([t]) = r(t)$. This gives $F_0(S_\Sigma)$ the stated $G$-operad structure, and the unit is $\pi^{-1} \circ \eta : S \to F_0(S)/\!\sim \,\, \to F_0(S_\Sigma)$.
\end{proof}

\section{The proof of \protect{Theorem \ref{thm:admfree}}}\label{subsec:pfkeycalc}

In this section, we perform the key calculation of \S\ref{sec:realize}. For readability, we begin by recalling some concepts and notation, and then we prove the following result.

\begin{thmkeycalc*} If $S$ is a  $N$ symmetric sequence, then $A(F(S)) = \la A (S) \ra$.
\end{thmkeycalc*}

\subsection{Recollections} Fix a finite group $G$.

If $T$ is a finite $H$-set, then the \emph{graph subgroup} $\Gamma(T) \subset G \times \Sigma_{\abs{T}}$ is the graph of some permutation representation of $T$ (Definition \ref{defn:extnorm}). The subgroup $\Gamma(T)$ is well-defined up to conjugation, and it is canonically determined if $T$ has an order.

A \emph{$N$ symmetric sequence} in $\b{Set}^G$ is a $\Sigma$-free symmetric sequence $X$ such that $X(0)^G , X(2)^G \neq \varnothing$. We say that \emph{$X$ admits $T$} if $X(\abs{T})^{\Gamma(T)} \neq \varnothing$, and we write $A(X)$ for the class of admissible sets of $X$ (Definition \ref{defn:Nsym}). The class $A(X)$ is a \emph{coefficient system}. This means it is closed under isomorphism, restriction, and conjugation by elements of $G$ (Definition \ref{defn:coeffsys}).

A \emph{$N$ operad} is an operad in $\b{Set}^G$, which is also a $N$ symmetric sequence (Definition \ref{defn:NOp}). Every $N$ symmetric sequence generates a free $N$ operad $F(S)$, and if $\s{O}$ is any $N$ operad, then $A(\s{O})$ is an indexing system. This is a coefficient system that contains all trivial actions, and is closed under subobjects, coproducts, and self-induction (Definition \ref{defn:indsys}). Every coefficient system $\c{C}$ generates an indexing system $\la \c{C} \ra$ (Definition \ref{defn:genind}).

Theorem \ref{thm:admfree} asserts that taking admissible sets commutes with free generation. This is a computation of the fixed points of a free $N$ operad $F(S)$. We shall see that the indexing system axioms mirror the structure of composition in $F(S)$.

The free operad $F(S)$ is typically defined as a large colimit, but it is hard to compute the fixed points of a quotient. Therefore we shall use a different model, denoted $F_0(S_\Sigma)$. This operad is described in detail in \S\ref{sec:quotfreeop}. We recommend rereading \S\ref{subsec:sumfree}, but briefly, $S_\Sigma$ is a set of $\Sigma$-orbit representatives for $S$, and the elements $t \in F_0(S_\Sigma)$ are formal composites of operations in $S_\Sigma$. The $\Sigma$-action permutes inputs, and the $G$-action is computed by conjugating every operation, and then replacing operations with their representatives in $S_\Sigma$ (Construction \ref{const:F0} and Theorem \ref{thm:freeGsigfree}). There is a related operad $F_0(S)$, whose elements are formal composites of operations in $S$. It has the same nonequivariant operad structure, but its $G$-action is just conjugation (Construction \ref{const:F0}).

Given any $t \in F_0(S_\Sigma)$, the \emph{complexity} of $t$ is the length of the longest chain of nested parentheses in $t$ (Definition \ref{defn:opterms}). Thus if $t = f(t_1,\dots,t_n)$, then the complexity of each $t_i$ is less than the complexity of $t$. We write $\ol{t}$ for $t$, but with all variables reindexed as $x_1,x_2,\dots$ from left to right (Notation \ref{nota:tbar}).

\subsection{The proof of Theorem \ref{thm:admfree}} As explained in \S\ref{subsec:keycalc}, the inclusion $A(F(S)) \supset \la A(S) \ra$ is easy. We now consider $A(F(S)) \subset \la A(S) \ra$.

\begin{defn} Suppose $t \in F_0(S_\Sigma)$, $H \subset G$ is a subgroup, and $T$ is a finite $H$-set. We say that $T$ is \emph{$t$-admissible} if $t \in F_0(S_\Sigma)(\abs{T})^{\Gamma(T)}$, where $\Gamma(T)$ is the graph of some permutation representation of $T$.
\end{defn}

We shall prove the following:
	\begin{equation*}\tag{$*$}
	\t{For any $t \in F_0(S_\Sigma)$ and finite $H$-set $T$, if $T$ is $t$-admissible, then $T \in \la A(S) \ra$.}
	\end{equation*}
Since every admissible set of $F_0(S_\Sigma)$ is $t$-admissible for some $t \in F_0(S_\Sigma)$, this will establish the inclusion $A(F(S)) = A(F_0(S_\Sigma)) \subset \la A(S) \ra$.

\begin{proof} We argue by induction on the complexity of $t \in F_0(S_\Sigma)$. If $t$ has complexity $0$, then $t = x_1 \in F_0(S_\Sigma)(1)$. Therefore every $t$-admissible set $T$ is an action of a subgroup $H \subset G$ on a point. It follows $T \in \la A(S) \ra$, because indexing systems contain all trivial actions.

Now suppose $t = f(t_1,\dots,t_n) \in  F_0(S_\Sigma)$ for some $f \in S_\Sigma (n)$ and $t_1,\dots,t_n$. Assume ($*$) is true for all $t$ of smaller complexity. For any $1 \leq i \leq n$, the complexity of $t_i$ is less than the complexity of $t$ and equal to the complexity of $\ol{t_i} \in F_0(S_\Sigma)$, so by induction, every $\ol{t_i}$-admissible set is contained in $\la A(S) \ra$.

Consider a $t$-admissible $H$-set $T$. We must prove that $T \in \la A(S) \ra$. The strategy is to use  the action on $F_0(S_\Sigma)$ to express $T$ in terms of $\ol{t_i}$-admissible sets. Since indexing systems are closed under isomorphism, we may assume that $T = \{1,\dots, \abs{T}\}$ and $t \in F_0(S_\Sigma)(\abs{T})^{\Gamma(T)}$, where $\Gamma(T) = \{(h,\sigma(h)) \, | \, h \in H\}$ and $\sigma(h) = h \cdot (-) : T \to T$. Therefore $h * t \cdot \sigma(h)^{-1} = (h,\sigma(h)) * t = t$, and hence
	\[
	h * t = t \cdot \sigma(h)	\quad\quad\t{(for all $h \in H$)}.
	\]
This is important. By Theorem \ref{thm:freeGsigfree}, the term $h * t$ is computed by multiplying in $F_0(S)$ and then shuffling subterms of $t$ around, whereas the term $t \cdot \sigma(h)$ is computed by permuting the variables of $t$ according to $\sigma(h) = h \cdot (-) : T \to T$. Thus, we can analyze the $H$-action on $T$ using the recursive definition of $h * (-)$.

For every $h \in H$, write $h \cdot f = f_h \cdot \tau(h)$ for unique $f_h \in S_\Sigma (n)$ and $\tau(h) \in \Sigma_n$. Here $f$ is the first letter of $t$, and products are computed in $F_0(S)$. Then
	\[
	f_h(h * t_{\tau(h)^{-1}1} , \dots , h* t_{\tau(h)^{-1}n}) = h * f(t_1,\dots,t_n) = f(t_1,\dots,t_n) \cdot \sigma(h).
	\]
The first letters must agree, so $f_h = f$ and $h \cdot f = f \cdot \tau(h)$. Hence $(h,\tau(h)) \cdot f = h \cdot f \cdot \tau(h)^{-1} = f$ for all $h \in H$, which implies the subgroup $\{ (h,\tau(h)) \, | \, h \in H \} \subset G \times \Sigma_n$ fixes $f \in S(n)$. Since $S$ is $\Sigma$-free, the set $\{ (h,\tau(h)) \, | \, h \in H \}$ is the graph subgroup $\Gamma(U)$ of an $H$-set $U$ with permutation representation $\tau : H \to \Sigma_n$. Thus $U \in A(S) \subset \la A(S) \ra$. Decomposing $U$ into orbits, we see that $H/K \in \la A(S) \ra$ for every suborbit $H/K \subset U$, because indexing systems are closed under subobjects.

Now we group the variables in each $t_i$ along the orbits of $U$. Let 
	\[
	T_i = \{ j \in \bb{N} \, | \, x_j \t{ appears in } t_i\}	\quad\quad ( 1 \leq i \leq n ),
	\]
so that $T_1 \sqcup \cdots \sqcup T_n = T$ as sets. For each orbit $O \subset U = \{1,\dots,n\}$, let
	\[
	T_O = \coprod_{i \in O} T_i.
	\]
We claim that $T_O$ is a sub-$H$-set of $T$. In fact, we shall show $\sigma(h)(T_i) = T_{\tau(h)i}$.

For any $h \in H$, write $t_i'$ for the term obtained from $t_i$ by replacing each variable $x_i$ with $x_{\sigma(h)^{-1}i}$. Then
	\[
	f(t_1',\dots,t_n') = t \cdot \sigma(h) = h * t = f(h * t_{\tau(h)^{-1}1},\dots, h * t_{\tau(h)^{-1}n}) ,
	\]
and therefore $t_{\tau(h)i}' = h * t_{i}$ by Proposition \ref{prop:parseterm}. Thus, the same variables appear in $t'_{\tau(h)i}$ and $t_i$, which means $\sigma(h)^{-1}(T_{\tau(h)i}) = T_i$. This proves that $T_O$ is a sub-$H$-set of $T$. Moreover, there is an isomorphism
	\[
	T \cong \coprod_{O} T_O
	\]
of $H$-sets. Thus, to prove $T \in \la A(S) \ra$, it will be enough to show $T_O \in \la A(S) \ra$ for each orbit $O$, because indexing systems are closed under coproducts.

Consider $T_{H/K} = \coprod_{aK \in H/K} T_{aK}$ for a given orbit $H/K \subset U$. Then $\sigma(h)(T_{eK}) = T_{hK}$ for each $h \in H$. Thus, $T_{eK}$ is a sub-$K$-set of $\t{res}^H_K T_{H/K}$ that generates $T_{H/K}$ as an $H$-set, and $\abs{T_{H/K}} = \abs{H:K} \cdot \abs{T_{eK}}$. Therefore the inclusion $T_{eK} \hookrightarrow \t{res}^H_K T_{H/K}$ induces an isomorphism
	\[
	T_{H/K} \cong \t{ind}_K^H T_{eK}.
	\]
Thus, to prove $T_{H/K} \in \la A(S) \ra$, it will be enough to show $T_{eK} \in \la A(S) \ra$, because $H/K \in \la A(S) \ra$ and indexing systems are closed under self-induction.

However, the $K$-action on $T_{eK}$ is isomorphic to the $K$-action on the variables of one of the subterms $t_i$ in $f(t_1,\dots,t_n)$, and this is isomorphic to the $K$-action on the variables of $\ol{t_i}$. This $K$-action is $\ol{t_i}$-admissible, by the definition of the $G \times \Sigma$-action on $F_0(S_\Sigma)$, and therefore $T_{eK} \in \la A(S) \ra$ by the induction hypothesis.

Thus $T_{eK}$, $T_{H/K} \cong \t{ind}_K^H T_{eK}$, and $T \cong \coprod_{O} T_O$ are all elements of $\la A(S) \ra$, which is what we needed to prove. By induction on the complexity of $t \in F_0(S_\Sigma)$, we conclude that $A(F(S)) \subset \la A(S) \ra$.
\end{proof}

\section{Equivariant Barratt-Eccles operads}\label{sec:assocNop}

In \S\ref{sec:realize}, we showed how to realize every indexing system as a free $N$ operad. In this section, we construct strictly associative and unital realizations (Definition \ref{defn:AsT}). These are the smallest models of $N_\infty$ operads that we know of, and after applying the functor $E : N\t{-}\b{Op}^G \to N_\infty\t{-}\b{Op}^G$, they become $N_\infty$ variants of the Barratt-Eccles operad. We summarize the basic properties of these operads in Theorem \ref{thm:assocNop}, and then we analyze their combinatorics in \S\ref{defn:AsT}. The proof of Theorem \ref{thm:assocNop} is given in \S\ref{subsec:pfassocNop}.

\subsection{Associative $N$ operads} For each indexing system $\c{I}$, we construct an associative and unital operad $\b{As}(\c{I})$ as follows.

\begin{defn}\label{defn:AsT} Let $\c{T} = (T_{\alpha})_{\alpha \in J}$ and $\b{F}_{\c{T}}$ be as in Definition \ref{defn:FT}, and suppose $\eta : S_{\c{T}} \to \b{F}_{\c{T}}$ is the unit of the adjunction. Write
	\[
	e = \eta(G \times \{\t{id}_0\}) \, , \, \otimes = \eta(G \times \{\t{id}_2\}) \, , \, \t{and } \begin{array}{c}\bigotimes_{T_{\alpha}}\end{array} \!\!\!\!\!\! = \eta(\Gamma(T_{\alpha}))
	\]
for every index $\alpha \in J$.  We define $\b{As}_{\c{T}}$ to be the quotient
	\[
	\b{As}_{\c{T}} = \frac{\b{F}_{\c{T}}}{
	\Bigg\la 
	\begin{array}{c}
	\gamma(\otimes;\otimes,\t{id}) \sim \gamma(\otimes;\t{id};\otimes) \, , \,
	\gamma(\otimes;e,\t{id}) \sim \t{id} \sim \gamma(\otimes ; \t{id} , e)	\\
	\gamma(\bigotimes_{T_{\alpha}} ; e,\dots,e) \sim e	 \, , \,
	\gamma(\bigotimes_{T_{\alpha}} ; e, \dots , e , \t{id} , e , \dots , e) \sim \t{id}
	\end{array}
	\Bigg| \,\,
	\alpha \in J
	\Bigg\ra 
	} 
	\]
of $\b{F}_{\c{T}}$ by the indicated relations. In $\gamma(\bigotimes_{T_{\alpha}} ; e, \dots , e , \t{id} , e , \dots , e)$, we allow $\t{id}$ to range over the $2$nd -- $[\abs{T}+1]$st arguments of $\gamma$. If $\abs{T_{\alpha}} = 0$, then we understand the lower left relation to be $\bigotimes_{T_{\alpha}} \sim e$. 

For any indexing system $\c{I}$, let $\b{As}(\c{I}) = \b{As}_{\b{O}(\c{I})}$, where $\b{O}(\c{I})$ is the set of all nontrivial orbits $H/K \in \c{I}$. 
\end{defn}

The operads $\b{As}(\c{I})$ have a number of useful properties, which are summarized in the theorem below.

\begin{thm}\label{thm:assocNop} The functor $A : N\t{-}\b{Op}^G \to \b{Ind}(G)$ has a functorial section
	\[
	\b{As} : \b{Ind}(G) \to N\t{-}\b{Op}^G
	\]
such that
	\begin{enumerate}
		\item{}$\b{As}(\ub{triv})$ is the associativity operad equipped with a trivial $G$-action, 
	\end{enumerate}
and for all $\c{I} \in \b{Ind}$:
	\begin{enumerate}[resume]
		\item{}$\b{As}(\c{I})$ is finitely presented, and
		\item{}$\abs{\b{As}(\c{I})(0)} = \abs{\b{As}(\c{I})(1)} = 1$, and there is $C = C(\c{I}) \in \bb{N}$ such that for every $n \geq 2$, we have the inequality $\abs{\b{As}(\c{I})(n)} < C^n (n!)^2$.
	\end{enumerate}
\end{thm}
The proof will be given in \S\ref{subsec:pfassocNop}. For now, we explain the significance of this result. Functoriality of $\b{As}$ in $\c{I}$ implies we can restrict $\b{As}(\c{J})$ actions to $\b{As}(\c{I})$ actions directly, provided that $\c{I} \subset \c{J}$. This eliminates the need to pass through a zig-zag $\b{As}(\c{I}) \stackrel{\sim}{\to} \b{As}(\c{I}) \times \b{As}(\c{J}) \to \b{As}(\c{J})$.

Condition (1) says that $\b{As}(\c{I})$ is a generalization of the usual associative operad.

Conditions (2) and (3) are bounds on the size of $\b{As}(\c{I})$, but first, a bit of context. Recall that the categorical Barratt-Eccles operad $\s{P}$ has $n$th space $\til{\b{As}(n)}$, where $\b{As}(n) = \Sigma_n$ is the associativity operad, and $\til{(-)} : \b{Set} \to \b{Cat}$ is the right adjoint to the object functor. 

In their work on equivariant infinite loop space theory, Guillou-May-Merling-Osorno consider the coinduced operad $\s{P}_G(n) \cong \til{\b{Set}(G,\Sigma_n)}$. This is a genuine $E_\infty$ $G$-operad, meaning it is $N_\infty$ and its indexing system is $\ub{Set}$. The operad $\s{P}_G$ was thought by many to be the smallest model for an $E_\infty$ $G$-operad, because $\s{P}$ is certainly the smallest model nonequivariantly.

This intuition is false. Work in \cite{BBKOOTX} shows that $\b{Set}(G,\b{As}) = \t{Ob}(\s{P}_G)$ is not finitely generated when $G$ is nontrivial, and if $\c{I} = \ub{Set}$, then (3) implies
	\[
	\underset{n \to \infty}{\t{lim}} \frac{\abs{\b{As}(\ub{Set})(n)}}{\abs{\t{Ob}(\s{P}_G(n))}} = 0
	\]
whenever $\abs{G} > 2$. The bound on $\abs{\b{As}(\c{I})(0)}$ is also useful. It says that $\b{As}(\c{I})$ is a reduced operad, and therefore $E\b{As}(\c{I})$ is, too. This can be quite convenient in applications (cf. \cite[Remark 2.7]{BH3}).

We round off this section by proposing two new definitions.

\begin{defn}Let $\c{I}$ be a $G$-indexing system. The \emph{$\c{I}$-permutativity operad} is
	\[
	\s{P}(\c{I}) = \til{\b{As}(\c{I})},
	\]
where $\til{(-)} : \b{Set} \to \b{Cat}$ is the right adjoint to $\t{Ob} : \b{Cat} \to \b{Set}$. The \emph{$\c{I}$-Barratt-Eccles operad} is
	\[
	\s{E}(\c{I}) = E\b{As}(\c{I}),
	\]
where $E = B \circ \til{(-)}$ is the composite of $\til{(-)}$ and the classifying space functor.
\end{defn}

\begin{rem} The operad $\s{P}_G$ is homotopy terminal, and early attempts at Conjecture \ref{conj:BH} sought to realize arbitrary indexing systems as suboperads of $\s{P}_G$. Bonventre proved this is impossible  \cite[Example B.2.1]{BonThesis}, and the construction of $N_\infty$ permutativity operads has been a sticking point ever since. Our operads $\s{P}(\c{I})$ are one possible solution.
\end{rem}

\subsection{Identifying associative $N$ operads}\label{subsec:assocNopcomb} The proof of Theorem \ref{thm:assocNop} requires a precise description of $\b{As}_{\c{T}}$. This section works out the details. We recommend skimming it on a first reading.

\begin{lem}\label{lem:AsTpres} For any indexed set $(T_\alpha)_{\alpha \in J}$ of finite $G$-subgroup actions, the operad $\b{As}_{\c{T}}$ is isomorphic to the operad
	\[
    	\frac{F \Big( \coprod_{n \geq 0} \frac{G \times \Sigma_n}{G \times \{\t{id}_n\}} \sqcup \coprod_{\alpha \in J} \frac{G \times \Sigma_{\abs{T_{\alpha}}}}{\Gamma(T_{\alpha})}  \Big)}
	{\Bigg\la \begin{array}{c}
		\gamma( \Pi_m ; \t{id} , \dots , \Pi_n , \dots , \t{id} ) \sim \Pi_{m+n-1} \,\, , \,\, \Pi_1 \sim \t{id} \\
		\gamma(\bigotimes_{T_{\alpha}} ; \Pi_0,\dots, \Pi_0) \sim \Pi_0 \,\, , \,\, \gamma(\bigotimes_{T_{\alpha}} ; \Pi_0 , \dots, \t{id} , \dots , \Pi_0) \sim \t{id}
	\end{array} 
	\Bigg| 
		\begin{array}{c}
			m \geq 1 \, , \, n \geq 0	\\
			\alpha \in J
		\end{array}
	\Bigg\ra
	}
	,
	\]
where $\Pi_k = G \times \{\t{id}_k\}$ for all $k \geq 0$, and $\bigotimes_{T_{\alpha}} = \Gamma(T_{\alpha})$ for all $\alpha \in J$. If $\abs{T_{\alpha}} = 0$, we understand the bottom left relation to be $\bigotimes_{T_{\alpha}} \sim \Pi_0$.
\end{lem}

\begin{proof} The inclusion of generators 
	\[
	\coprod_{n = 0,2} \frac{G \times \Sigma_n}{G \times \{ \t{id}_n \}} \sqcup \coprod_{\alpha \in J} \frac{G \times \Sigma_{\abs{T_{\alpha}}}}{\Gamma(T_{\alpha})} 
	\,\,\hookrightarrow\,\,
	 \coprod_{n \geq 0} \frac{G \times \Sigma_n}{G \times \{\t{id}_n\}} \sqcup \coprod_{\alpha \in J} \frac{G \times \Sigma_{\abs{T_{\alpha}}}}{\Gamma(T_{\alpha})}
	\]
induces an isomorphism.
\end{proof}

The presentation of $\b{As}_{\c{T}}$ in Lemma \ref{lem:AsTpres} is easier to work with, because the relations are clearly ``reductions.'' We use it to solve the word problem for $\b{As}_{\c{T}}$.

\begin{prop}\label{prop:AsTrep} Let $\c{T} = (T_{\alpha})_{\alpha \in J}$ be an indexed set of finite $G$-subgroup actions. The operad $\b{As}_{\c{T}}$ is isomorphic to a sub-symmetric sequence of the free operad
	\[
	\s{F}_{\c{T}} = F \Bigg( \coprod_{n \geq 0} \frac{G \times \Sigma_n}{G \times \{\t{id}_n\}} \sqcup \coprod_{\alpha \in J} \frac{G \times \Sigma_{\abs{T_{\alpha}}}}{\Gamma(T_{\alpha})}  \Bigg) ,
	\]
equipped with a reduced composition operation.
\end{prop}

\begin{proof} For each subgroup $H \subset G$, choose a set $\{e = r^H_1,\dots,r^H_{\abs{G:H}} \}$ of $G/H$-coset representatives once and for all. Then
	\[
	P_\Sigma = \coprod_{n \geq 0} \Bigg\{ G \times \{\t{id}_n\} \Bigg\} \sqcup \coprod_{\alpha \in J} \Bigg\{ r^H_i \Gamma(T_{\alpha}) 
	\,\, \Bigg|  
	\begin{array}{c}
		H \subset G , \, T_{\alpha} \t{ an $H$-set},	\\
		\t{and } 1 \leq i \leq \abs{G:H}
	\end{array}
	\Bigg\}
	\]
is a set of $\Sigma$-orbit representatives for the generators of $\s{F}_{\c{T}}$. It follows from Theorem \ref{thm:freeGsigfree} that $\s{F}_{\c{T}} \cong F_0(P_\Sigma)$ with a twisted $G$-action.

We identify the congruence relation $\sim$ on $F_0(P_\Sigma)$ that is generated by the relations in Lemma \ref{lem:AsTpres}. For any $n \geq 0$ and $t,t' \in F_0(P_\Sigma)(n)$, declare $t \to t'$ if $t'$ is obtained by replacing a subterm $s$ of $t$ with a new subterm $s'$, in one of the following ways:
	\[
	\begin{array}{|c|c|}
	\hline
	s	&	s'	\\
	\hline
	\Pi_m(t_1,\dots, t_{i-1},\Pi_n(t_i,\dots,t_{i+n-1}),t_{i+n},\dots,t_{m+n-1})	&	\Pi_{m+n-1}(t_1,\dots,t_{m+n-1})	\\
	\Pi_1(t_1)	&	t_1	\\
	r^H_i\bigotimes_{T_{\alpha}}(\Pi_0(),\dots,\Pi_0())	&	\Pi_0()	\\
	r^H_i\bigotimes_{T_{\alpha}}(\Pi_0(),\dots,\Pi_0(),t_1,\Pi_0(),\dots,\Pi_0())	&	t_1	\\
	\hline
	\end{array}
	\]
In the first line, we require $m \geq 1$ and $n \geq 0$, and in the third and fourth lines, we require $\alpha \in J$ and $r^H_i$ to be a coset representative for $G/H$, where $H$ is the subgroup acting on $T_{\alpha}$. We say that $t$ is \emph{reduced} if there is no $t'$ such that $t \to t'$, and we write $rF_0(P_\Sigma) \subset F_0(P_\Sigma)$ for the sub-symmetric sequence of reduced elements.

Each of the substitutions above decreases the number
	\[
	w(t) = \#(\Pi_k \t{ symbols in $t$}) + 2 \cdot \#( r^i_H \!\!\! \begin{array}{c}\bigotimes_{T_\alpha}\end{array}\!\!\! \t{ symbols in $t$}) \geq 0 ,
	\]
and therefore $\to$ is well-founded. Write $\stackrel{*}{\to}$ for the reflexive and transitive closure of $\to$ and suppose $t \in F_0(P_\Sigma)$. By the diamond lemma, there is a unique $r(t) \in rF_0(P_\Sigma)$ such that $t \stackrel{*}{\to} r(t)$.
	
By Proposition \ref{prop:gencong}, the congruence relation $\sim$ is the equivalence relation generated by $\to$. It follows that $t \sim t'$ if and only if $r(t) = r(t')$, and therefore $rF_0(P_\Sigma)$ is a set of representatives for $\sim$, with $r(t)$ representing $t$. Hence
	\[
	\b{As}_{\c{T}} \cong F_0(P_\Sigma)/\!\!\sim \,\, \cong rF_0(P_\Sigma)
	\]
as symmetric sequences. Composition in $\b{As}_{\c{T}}$ is identified with $r \circ \gamma$, where $\gamma$ denotes composition in $F_0(P_\Sigma)$.
\end{proof}

Now we can estimate the size of $\b{As}_{\c{T}}$. We focus on $\c{T} = \b{O}(\c{I})$ for simplicity, but the same reasoning applies for any finite $\c{T}$.

\begin{lem}\label{lem:AsTcard} Suppose $\c{I}$ is an indexing system, $\b{O}(\c{I})$ is the set of nontrivial orbits in $\c{I}$, and write $\b{As}(\c{I}) = \b{As}_{\b{O}(\c{I})}$. Then $\abs{\b{As}(\c{I})(n)} = 1$ for $n = 0,1$, and there is a constant $C = C(\c{I}) \in \bb{N}$ such that $\abs{\b{As}(\c{I})(n)} < C^n (n!)^2$ for $n \geq 2$.
\end{lem}

\begin{proof}Keep notation as in the proof of Proposition \ref{prop:AsTrep} and set $\c{T} = \b{O}(\c{I})$. We count the number of elements in $rF_0(P_\Sigma)(n) \cong \b{As}(\c{I})(n)$ for each $n \geq 0$. The estimates are clear (and poor) when $\c{I} = \ub{triv}$, so assume $\c{I}$ is nontrivial.

Given $t \in F_0(P_\Sigma)(0)$, we can use the relation $r^H_i \bigotimes_{T_{\alpha}}() \sim \Pi_0()$ for empty $T_\alpha$'s to convert all nullary function symbols in $t$ into $\Pi_0$'s. Call the result $t'$. Now we use $\Pi_m(t_1,\dots,\Pi_0(),\dots,t_{m-1}) \sim \Pi_{m-1}(t_1,\dots,t_{m-1})$ and $r^H_i \bigotimes_{T_{\alpha}}(\Pi_0(),\dots,\Pi_0()) \sim \Pi_0()$ inductively to collapse $t'$ to $\Pi_0()$. Therefore $\abs{rF_0(P_\Sigma)(0)} = 1$.

The case for $rF_0(P_\Sigma)(1)$ is similar. We claim that every $t \in F_0(P_\Sigma)(1)$ can be reduced to $x_1$. For, if $t = f(t_1,\dots,t_n) = \gamma(\eta_0(f);\ol{t_1},\dots,\ol{t_n}) \cdot \sigma$, there is $1 \leq i \leq n$ such that $\ol{t_i}$ is unary and and $\ol{t_j}$ is nullary for $j \neq i$. By the above, we have $\ol{t_j} \sim \Pi_0()$, and we can assume $\ol{t_i} \sim x_1$ by induction on complexity. Therefore $t \sim f(\Pi_0(),\dots,x_1,\dots,\Pi_0()) \sim x_1$.

Now we make the estimate for $n \geq 2$. Every $t \in rF_0(P_\Sigma)(n)$ can be factored as $t = (b_m \circ_{k_{m-1}} b_{m-1} \circ_{k_{m-2}} \cdots \circ_{k_1} b_1) \cdot \sigma$, where $\sigma \in \Sigma_n$, $\circ_k$ is partial composition, and $b_1,\dots,b_m$ are \emph{basic terms} of the form
	\begin{align*}
		\Pi_2(x_1,x_2)	\quad\t{or}\quad 
		\!\!\!\begin{array}{c} r^H_i\bigotimes_{H/K}\end{array}\!\!\!(t_1,\dots,t_{\abs{H:K}}) ,
	\end{align*}
such that all of the terms $t_j$ are either variables or $\Pi_0()$'s, and at least two of the $t_j$'s are variables. The arity of each basic term is at least $2$. Hence
	\[
	2 \leq \abs{b_1} < \abs{b_2 \circ_{k_1} b_1} < \cdots < \abs{b_m \circ_{k_{m-1}} b_{m-1} \circ_{k_{m-2}} \cdots \circ_{k_1} b_1} =n ,
	\]
and it follows $m < n$. 

Let $B$ be the set of all basic terms and set $C = \abs{B} \geq 2$. For each $m = 1, \dots , n-1$, there are no more than $C^m$ choices of basic operations $(b_1,\dots,b_m)$ such that $\abs{b_1} + \cdots + \abs{b_m} + m - 1 = n$, and for each choice $(b_1,\dots,b_m)$, there are no more than $n!$ choices of sequences $(k_1,\dots,k_{m-1})$ such that $1 \leq k_j \leq  \abs{b_1} + \cdots + \abs{b_j} - j + 1$. Summing over $m$ and choosing a permutation $\sigma \in \Sigma_n$ shows there are fewer than $C^n (n!)^2$ $n$-ary expressions of the form $(b_m \circ_{k_{m-1}} b_{m-1} \circ_{k_{m-2}} \cdots \circ_{k_1} b_1) \cdot \sigma$.
\end{proof}

\subsection{The proof of Theorem \ref{thm:assocNop}}\label{subsec:pfassocNop} In this section, we prove Theorem \ref{thm:assocNop}, starting with a calculation of the admissible sets of $\b{As}_{\c{T}}$.

\begin{lem}\label{lem:admAs} For any $\c{T}$, $\b{As}_{\c{T}}$ is a $N$ operad, and $A(\b{As}_{\c{T}}) = \la T_{\alpha} \, | \, \alpha \in J \ra$.
\end{lem}

\begin{proof} Let $\s{F}_{\c{T}}$ be as in Proposition \ref{prop:AsTrep}. There is an embedding of symmetric sequences $\b{As}_{\c{T}} \hookrightarrow \s{F}_{\c{T}}$. Therefore $\b{As}_{\c{T}}$ is $\Sigma$-free and $A(\b{As}_{\c{T}}) \subset A(\s{F}_{\c{T}})$. On the other hand, Lemma \ref{lem:AsTpres} implies there is a quotient operad map $\s{F}_{\c{T}} \to \b{As}_{\c{T}}$. Therefore $\b{As}_{\c{T}}(n)^G \neq \varnothing$ and $A(\s{F}_{\c{T}}) \subset A(\b{As}_{\c{T}})$. This proves that $\b{As}_{\c{T}}$ is a $N$ operad, and 
	\[
	A(\b{As}_{\c{T}}) = A(\s{F}_\c{T}) = \la T_\alpha \, | \, \alpha \in J \ra
	\]
by Theorem \ref{thm:admfree}.
\end{proof}

Now we can prove the theorem.

\begin{proof}[Proof of Theorem \ref{thm:assocNop}] Define $\b{As}(-) : \b{Ind}(G) \to N\t{-}\b{Op}^G$ by $\b{As}(\c{I}) = \b{As}_{\b{O}(\c{I})}$, where $\b{O}(\c{I})$ is the set of nontrivial orbits $H/K \in \c{I}$. The same argument given in the proof of Theorem \ref{thm:freerealize} shows that $\b{As}(-)$ is functorial, and Lemma \ref{lem:admAs} shows
	\[
	A(\b{As}(\c{I})) = \la H/K \, | \, H/K \in \c{I} \t{ is nontrivial} \ra = \c{I}.
	\]
Therefore $\b{As}$ is a section of $A : N\t{-}\b{Op}^G \to \b{Ind}(G)$.

We have $\b{As}(\ub{triv}) = \b{As}$ by inspection, and $\b{As}(\c{I}) = \b{As}_{\b{O}(\c{I})}$ is finitely generated because $\b{O}(\c{I})$ is finite. Lemma \ref{lem:AsTcard} gives the desired cardinality bound.
\end{proof}

\section{Model categories of discrete $G$-operads}\label{sec:modstr}

This final section interprets \S\ref{subsec:Nop} and \S\ref{sec:realize} through a model categorical lens. We set up the basic model structures in \S\S\ref{susbec:modstrOpG}--\ref{subsec:LHNOp}, and then we compare our work to \cite{GutWhite} and \cite{BonPer} is \S\ref{subsec:compare}.

We have a few reasons for introducing this formalism. To start, we find it clarifying. The free operads in \S\ref{sec:realize} may seem ad hoc, but they are completely natural from a model categorical perspective (cf. Proposition \ref{prop:freerealcof}). Model categorical language also helps explain the relationship between our construction of $\b{F}(\c{I})$, and the realizations in \cite{GutWhite} and \cite{BonPer} (cf. \S\ref{subsec:compare}). That being said, the associative $N$ operads considered in \S\ref{sec:assocNop} do not mesh well with model structures. The operad $\b{As}_{\c{T}}$ is just too small to be cofibrant, and should be understood on the point-set level.

\begin{rem} Looking ahead, we will truly need these model structures in subsequent work. We could do things by hand in this paper, but parts of \cite{RubCateg} require a more sophisticated approach.
\end{rem}

\subsection{Model category structures on $\b{Op}^G$}\label{susbec:modstrOpG} A model category must be bicomplete, which implies we cannot literally equip the category $N\t{-}\b{Op}^G$ of $N$ operads with a model structure. Instead, we consider the category $\b{Op}^G = \b{Op}(\b{Set}^G)$ of all operads in $G$-sets, and then we cut things down later.

We start on the point-set level. The following holds in general (cf. \cite[\S2.3]{Rezk}).

\begin{lem} The category $\b{Op}^G$ is complete and cocomplete.
\end{lem}

Limits are computed levelwise in $\b{Set}^G$, and colimits are similar to colimits of nonabelian groups. We write $*$ for the coproduct in $\b{Op}^G$.

The category $\b{Op}^G$ also has a small set of small generators.

\begin{lem} The category $\b{Op}^G$ is locally finitely presentable.
\end{lem}

\begin{proof}The free operads $F(G \times \Sigma_n)$ form a strong generator of $\b{Op}^G$ \cite[\S0.6]{AR}, where $n \geq 0$ is a nonnegative integer. Moreover, each of the operads $F(G \times \Sigma_n)$ is finitely presentable. Therefore $\b{Op}^G$ is locally finitely presentable by \cite[Theorem 1.11]{AR}.
\end{proof}

Our ultimate goal is to construct a simplicial model category. We therefore give $\b{Op}^G$ a simplicial enrichment. There is an adjunction
	\[
	(-)_0 : \b{sSet} \rightleftarrows \b{Set} : E = N \circ \til{(-)}
	\]
where $(-)_0$ is the $0$-simplices functor, $\til{(-)} : \b{Set} \to \b{Cat}$ is the right adjoint to the object functor, and $N : \b{Cat} \to \b{sSet}$ is the nerve functor. As in \S\ref{subsec:Nop}, $E(\varnothing) = \varnothing$ and $EX \simeq *$ if $X \neq \varnothing$.

Since $(-)_0$ and $E$ are both limit-preserving functors, we may use the adjunction $(-)_0 \dashv E$ to enrich, tensor, and cotensor $\b{Op}^G$ over $\b{sSet}$ (cf. \cite[Theorem 3.7.11]{Riehl}).

\begin{lem}\label{lem:sSetenr} The category $\b{Op}^G$ is enriched, tensored, and cotensored over the category $\b{sSet}$ of simplicial sets, with:
	\begin{enumerate}[label=(\alph*)]
		\item{}hom objects: $\ub{O}\b{p}^G(\s{O},\s{O}') = E\b{Op}^G(\s{O},\s{O}')$,
		\item{}tensors: $K \otimes \s{O} = K_0 \cdot \s{O}$, the $K_0$-fold coproduct of copies of $\s{O}$, and
		\item{}cotensors: $\s{O}^K = \s{O}^{K_0}$, the $K_0$-fold product of copies of $\s{O}$,
	\end{enumerate}
where $\s{O},\s{O}' \in \b{Op}^G$ and $K \in \b{sSet}$.
\end{lem}

We could have done the same thing with $\b{Op}^G$ replaced by almost any $1$-category, but it is a reasonable choice for $\b{Op}^G$ because we are really thinking of $\s{O} \in \b{Op}^G$ as the categorical operad $\til{\s{O}}$. The hom object between $\til{\s{O}_1}$ and $\til{\s{O}_2}$ is naturally a $1$-category that is isomorphic to $\til{\b{Op}^G}(\s{O}_1,\s{O}_2)$.

Now we make $\b{Op}^G$ into a model category.

\begin{defn}Let $\c{I}$ be an indexing system and $\Gamma \subset G \times \Sigma_n$. We say that $\Gamma = \{ (h,\sigma(h)) \, | \, h \in H\}$ is an \emph{$\c{I}$-graph subgroup} if $\sigma : H \to \Sigma_n$ is the permutation representation of a member of $\c{I}$. A morphism $f : \s{O}_1 \to \s{O}_2$ in $\b{Op}^G$ is an \emph{$\c{I}$-weak equivalence} if $E{f} : E{\s{O}_1}(n)^{\Gamma} \to E{\s{O}_2}(n)^{\Gamma}$ is a weak homotopy equivalence of topological spaces for every $n \geq 0$ and $\c{I}$-graph subgroup $\Gamma \subset G \times \Sigma_n$.
\end{defn}

This boils down to the condition that $\s{O}_1(n)^{\Gamma}$ is nonempty whenever $\s{O}_2(n)^{\Gamma}$ is nonempty, provided that $\Gamma$ is an $\c{I}$-graph subgroup.

\begin{prop}\label{prop:Imodstr}Fix an indexing system $\c{I}$. The category $\b{Op}^G$, together with the $\c{I}$-weak equivalences, can be enhanced to a right proper, combinatorial, simplicial model category with generating cofibrations
	\[
	\s{I}_{\c{I}} = \Bigg\{ \{\t{id}\} \longrightarrow F \Bigg( \frac{G \times \Sigma_{n}}{\Gamma} \Bigg) \,\, \Bigg| \,\, 
	\begin{array}{c}
    		n \geq 0 , \, \Gamma \subset G \times \Sigma_n \t{ an}	\\
		\t{$\c{I}$-graph subgroup}
	\end{array}
	\Bigg\}
	\]
and generating acyclic cofibrations
	\[
	\s{J}_{\c{I}} = \Bigg\{ F \Bigg( \frac{G \times \Sigma_{n}}{\Gamma} \Bigg) \stackrel{i_0}{\longrightarrow} \Delta^1 \otimes F \Bigg( \frac{G \times \Sigma_{n}}{\Gamma} \Bigg)	
	\,\, \Bigg| \,\,
	\begin{array}{c}
    		n \geq 0 , \, \Gamma \subset G \times \Sigma_n \t{ an}	\\
		\t{$\c{I}$-graph subgroup}
	\end{array}
	\Bigg\} .
	\]
Here $F : (\b{Set}^G)^\Sigma \rightleftarrows \b{Op}^G : U$ is the free-forgetful adjunction, and $\{\t{id}\} \cong F(\varnothing)$ is the initial operad. Moreover: 
	\begin{enumerate}
		\item{}every object of $\b{Op}^G$ is $\c{I}$-fibrant, and
		\item{}every simplicial mapping space in $\b{Op}^G$ is either empty or contractible.
	\end{enumerate}
\end{prop}

\begin{proof} The construction of this model structure is a straightforward application of the small object argument (cf. \cite[Theorem 15.2.3]{MayPon}). It is also straightforward to verify that it is right proper and that every object is fibrant. The only interesting point is that axiom SM7 holds, which we now prove.

Suppose $i : \s{A} \to \s{X}$ is an $\c{I}$-cofibration and $p : \s{E} \to \s{B}$ is an $\c{I}$-fibration, and consider the map
	\[
	\ul{(i^*,p_*)} : \ub{O}\b{p}^G(\s{X},\s{E}) \to \ub{O}\b{p}^G(\s{A},\s{E}) \times_{\ub{O}\b{p}^G(\s{A},\s{B})} \ub{O}\b{p}^G(\s{X},\s{B}).
	\]
If either $i$ or $p$ is an $\c{I}$-weak equivalence, then $\ul{(i^* , p_*)}$ is a weak equivalence. Indeed, the domain and codomain are either empty or contractible, and if the codomain is nonempty, then the domain is nonempty by lifting. Thus, axiom SM7 will follow if we show that $\ul{(i^*,p_*)}$ is a Kan fibration.

By the adjunction $(-)_0 \dashv E$, the simplicial map $\ul{(i^* , p_*)}$ is a Kan fibration if and only if the set map $(i^*,p_*)$ has the right lifting property with respect to the inclusion $\{0\} \to \{0,1\}$. This is easy to check when $p$ is an $\c{I}$-fibration and $i$ is a relative $\s{I}_{\c{I}}$-cell complex $i_1 : \s{O} \to \s{O} * F(S)$. Passing to retracts proves the result for general $\c{I}$-cofibrations. Therefore $\b{Op}^G$ is a simplicial model category.
\end{proof}

We do not know if the $\c{I}$-model structure on $\b{Op}^G$ is left proper, because we do not know how to compute the fixed points of the relevant pushouts.

\begin{rem} There are analogous $\c{I}$-model structures on $\b{Op}(\b{sSet}^G)$ and $\b{Op}(\b{Top}^G)$ by the work in \cite{GutWhite} and \cite{BonPer}. The adjunction $(-)_0 : \b{sSet} \leftrightarrows \b{Set} : E$ induces a Quillen adjunction between the $\c{I}$-model structures on $\b{Op}(\b{sSet}^G)$ and $\b{Op}^G$ because $(-)_0$ sends generating (acyclic) cofibrations to (acyclic) cofibrations. In fact, one can construct the $\c{I}$-model structure on $\b{Op}^G$ by transport along $(-)_0 \dashv E$.
\end{rem}

\subsection{The homotopy theory of $N$ operads}\label{subsec:LHNOp} The $\ub{Set}$-model structure on $\b{Op}^G$ governs a broader homotopy theory than the homotopy theory of $N$ operads. One can prove that every bifibrant operad $\s{O} \in \b{Op}^G$ is $\Sigma$-free, but nothing ensures that $\s{O}(n)^G \neq \varnothing$. We fix things by passing to a slice category of $\b{Op}^G$. 

\begin{defn} Let $\b{F}$ be the free operad on $(G \times \Sigma_0)/G \sqcup (G \times \Sigma_2)/G$, and write $\b{Op}_{+}^G$ for the slice category $\b{F}/\b{Op}^G$ of symmetric operads in $\b{Set}^G$ under $\b{F}$.
\end{defn}

By adjunction, an object of $\b{Op}^G_{+}$ is the same thing as an operad $\s{O} \in \b{Op}^G$, equipped with marked operations $u \in \s{O}(0)^G$ and $p \in \s{O}(2)^G$. A morphism in $\b{Op}^G_{+}$ is just a morphism in $\b{Op}^G$ that preserves the markings.

We enrich, tensor, and cotensor $\b{Op}_{+}^G$ over $\b{sSet}$ as before, i.e. we declare $\ub{O}\b{p}^G_+(\s{O}_1,\s{O}_2) = E\b{Op}^G_+(\s{O}_1,\s{O}_2)$ and we define tensors and cotensors by adjunction (cf. Lemma \ref{lem:sSetenr}). From here, we use the $\ub{Set}$-model structure on $\b{Op}^G$ to create a model structure on $\b{Op}_{+}^G$. We summarize its properties.

\begin{thm}\label{thm:modOpG+} The category $\b{Op}^G_+$ is a right proper, combinatorial, simplicial model category. A morphism $f : \s{O}_1 \to \s{O}_2$ in $\b{Op}^G_+$ is a weak equivalence, fibration, or cofibration if, after forgetting markings, it is such a map in the $\ub{Set}$-model structure on $\b{Op}^G$. The generating cofibrations and acyclic cofibrations of $\b{Op}^G_+$ are the sets $\b{F} * \s{I}_{\ub{Set}}$ and $\b{F} * \s{J}_{\ub{Set}}$, where $\s{I}_{\ub{Set}}$ and $\s{J}_{\ub{Set}}$ are the corresponding generators for $\b{Op}^G$. Moreover:
	\begin{enumerate}
		\item{}every object of $\b{Op}^G_+$ is fibrant, 
		\item{}every cofibrant object of $\b{Op}^G_+$ is a $N$ operad (but not conversely), and 
		\item{}every mapping space in $\b{Op}^G_+$ is either empty or contractible.
	\end{enumerate}
\end{thm}

\begin{proof} The $\ub{Set}$-model structure on $\b{Op}^G$ lifts to a model structure on $\b{Op}^G_+ = \b{F}/\b{Op}^G$ by \cite[Theorem 15.3.6]{MayPon}, and the remaining claims about the unenriched model structure are standard. Axiom SM7 holds for $\b{Op}^G_+$, because for any cofibration $i : \s{A} \to \s{X}$ and fibration $p : \s{E} \to \s{B}$ in $\b{Op}^G_+$, the map	
	\[
	\ul{(i^*,p_*)} : \ub{O}\b{p}^G_+(\s{X},\s{E}) \to \ub{O}\b{p}^G_+(\s{A},\s{E}) \times_{\ub{O}\b{p}^G_+(\s{A},\s{B})} \ub{O}\b{p}^G_+(\s{X},\s{B})
	\]
is a pullback of the analogous map for $\ub{O}\b{p}^G$. It remains to show that every cofibrant operad $\s{O} \in \b{Op}^G_+$ is a $N$ operad.

If $\s{O} \in \b{Op}^G_+$ is cofibrant, then $\b{F} \hookrightarrow U\s{O}$ is a $\ub{Set}$-cofibration in $\b{Op}^G$, and since $\b{F}$ is $\ub{Set}$-cofibrant, so too is $U\s{O}$. Therefore $U\s{O}$ is a retract of a free operad $F(S)$ on a $\Sigma$-free symmetric sequence $S$. By universality, $F(S)$ must be $\Sigma$-free, and therefore $U\s{O}$ is also $\Sigma$-free because there is a map $U\s{O} \to F(S)$. It follows that $U\s{O}$ is a $N$ operad because we have another map $\b{F} \to U\s{O}$.
\end{proof}

Part (2) of Theorem \ref{thm:modOpG+} lets us relate $\b{Op}^G_+$ to $N\t{-}\b{Op}^G$.

\begin{prop}\label{prop:LHNOp} The cofibrant replacement functor $Q : \b{Op}^G_+ \to N\t{-}\b{Op}^G$ induces a Dwyer-Kan equivalence between the hammock localizations of $\b{Op}^G_+$ and $N\t{-}\b{Op}^G$. Therefore the functor $\bb{L}E = E \circ Q : \b{Op}^G_+ \to N_\infty\t{-}\b{Op}^G$ also induces a Dwyer-Kan equivalence between the corresponding hammock localizations.
\end{prop}

\begin{proof} Consider the functors below.
	\[
	\begin{tikzpicture}
		\node(1) at (0,0) {$N\t{-}\b{Op}^G$};
		\node(2) at (3,0) {$N\t{-}\b{Op}^G_{free}$};
		\node(3) at (6,0) {$(\b{Op}^G_+)_{cell}$};
		\node(4) at (9,0) {$\b{Op}^G_+$};
		\path[->]
		(1.9) edge [above] node {$F$} (2.172)
		(2.187) edge [below] node {$i$} (1.-9)
		(2.7) edge [above] node {$\b{F}*(-)$} (3.173)
		(3.187) edge [below] node {$U$} (2.-7)
		(3.7) edge [above] node {$i$} (4.168)
		(4.192) edge [below] node {$Q$} (3.-7)
		;
	\end{tikzpicture}
	\]
Here $(N\t{-}\b{Op})_{free}$ is the full subcategory of $\b{Op}^G$ spanned by frees, $(\b{Op}^G_+)_{cell}$ is the full subcategory spanned by cell complexes, $i$ denotes inclusion, $U$ is forgetful, $F$ is free, and $Q$ is cofibrant replacement. Every composite of opposing pairs is naturally weakly equivalent to the identity. Therefore all six of these functors induce Dwyer-Kan equivalences by \cite[\S 3]{DK}. The same is true for $\bb{L}E$ by Theorem \ref{thm:Gsetmod}.
\end{proof}

Since every mapping space in $L^H(\b{Op}^G_+)$ is empty or contractible, we deduce the same holds for $N\t{-}\b{Op}^G$ and $N_\infty\t{-}\b{Op}^G$.

\begin{cor}\label{cor:LHNop} Every mapping space in the hammock localization $L^H(N\t{-}\b{Op}^G)$ is either empty or contractible, and the same is true for $L^H(N_\infty\t{-}\b{Op}^G)$.
\end{cor}

This reproves \cite[Proposition 5.5]{BH1}. We end this section with an observation.

\begin{rem} Consider the functor $A : \t{Ho}(N\t{-}\b{Op}^G) \to \b{Ind}(G)$ once more. Corollary \ref{cor:LHNop} implies that $A$ is faithful, and Theorems \ref{thm:freerealize} and \ref{thm:assocNop} imply that $A$ is surjective. Fullness can deduced be using the product trick. If $\s{O}_1$ and $\s{O}_2$ are $N$ operads and $A(\s{O}_1) \subset A(\s{O}_2)$, then $\s{O}_1 \stackrel{\sim}{\leftarrow} \s{O}_1 \times \s{O}_2 \to \s{O}_2$ represents a morphism in $\t{Ho}(N\t{-}\b{Op}^G)$ that lifts the inclusion. This is a purely combinatorial proof that $A : \t{Ho}(N\t{-}\b{Op}^G) \to \b{Ind}(G)$ is an equivalence. Thus, the only topological ingredient in the our proof of the classification of $N_\infty$ operads (Theorem \ref{thm:classNinfty}) is the equivalence between $N_\infty$ operads and $N$ operads (Theorem \ref{thm:Gsetmod}).
\end{rem}

\subsection{Comparisons of $N_\infty$ realizations}\label{subsec:compare}

In \S\ref{sec:realize}, we showed how to realize arbitrary indexing systems using the free $N$ operads $\b{F}_{\c{T}}$. We now explain how to compare these operads to the operads constructed in \cite{GutWhite} and \cite{BonPer}. Recall that $\b{Com}$ is the terminal operad, whose levels are $\b{Com}(n) = *$ for all $n \geq 0$.

\begin{prop}\label{prop:freerealcof}The $N$ operads $\b{F}_{\b{O}(\c{I})}$ and $\b{F}_{\b{N}(\c{I})}$, described in Theorem \ref{thm:freerealize} and Example \ref{ex:choosegens}, are cofibrant replacements of the operad $\b{Com}$ in the $\c{I}$-model structure on $\b{Op}^G$.
\end{prop}

\begin{proof} Let $\s{F} = \b{F}_{\b{O}(\c{I})}$ or $\b{F}_{\b{N}(\c{I})}$. The operad $F(G \times \Sigma_n/\Gamma)$ is $\c{I}$-cofibrant for every $\c{I}$-graph subgroup $\Gamma$, and $\s{F}$ is a coproduct of such operads. Therefore $\s{F}$ is also $\c{I}$-cofibrant. Moreover, the unique morphism $\s{F} \to \b{Com}$ is an $\c{I}$-acyclic fibration, because Theorem \ref{thm:admfree} enusres $A(\s{F}) = \c{I}$.
\end{proof}

Thus, the functor $\b{F} : \b{Ind} \to N\t{-}\b{Op}$ in Theorem \ref{thm:freerealize} constructs operads that are formally analogous to Guti\'{e}rrez and White's $N_\infty$ operads \cite[Theorem 4.7]{GutWhite}. They prove that an $\c{I}$-cofibrant replacement of the operad $\b{Com} \in \b{Op}(\b{Top}^G)$ is a $N_\infty$ realization of $\c{I}$.

More concretely, consider the $N_\infty$ operad $E{\b{F}}_{\b{N}(\c{I})}$. It is constructed by generating a free, discrete operad $\b{F}_{\b{N}(\c{I})}$ with all operations specified by $\c{I}$, and then killing all homotopy groups with $E$. Guti\'{e}rrez and White's operads are similarly constructed. By the small object argument, an $\c{I}$-cofibrant replacement of $\b{Com}$ may be presented as a transfinite sequential colimit $\s{O}_{\c{I}} = \t{colim}_{\alpha < \gamma} \s{O}_{\alpha}$, where
	\begin{enumerate}[label=(\roman*)]
		\item{}$\s{O}_0 = \{\t{id}\}$, 
		\item{}$\s{O}_{\alpha + 1}$ is obtained from $\s{O}_{\alpha}$ by attaching a free cell $F( (G \times \Sigma_{n}/\Gamma) \times D^m)$ along every operad map $F( (G \times \Sigma_{n}/\Gamma) \times S^{m-1}) \to \s{O}_\alpha$, where $m , n \geq 0$ and $\Gamma$ is an $\c{I}$-graph subgroup, and		\item{}$\s{O}_{\beta} = \t{colim}_{\alpha < \beta} \s{O}_{\alpha}$ for each limit ordinal $\beta < \gamma$.
	\end{enumerate}
In particular, $\s{O}_1$ splits as $F(\coprod_{\Gamma} G \times \Sigma_n/\Gamma) * \s{O}_1'$, where $\Gamma$ ranges over all $\c{I}$-graph subgroups and $\s{O}_1'$ is built from $F((G \times \Sigma_1)/H \times D^m)$-cell attachments. Subsequent stages introduce more generators and kill elements of homotopy. By compactness, all homotopy is killed in the limit.

Bonventre and Pereira \cite[Remark 6.73]{BonPer} also construct $N_\infty$ operads as cofibrant replacements of $\b{Com}$, but they use a different model. Their powerful theory realizes the indexing system $\c{I}$ as a monadic bar construction $\s{B}_{\c{I}} = B_\bullet (\widehat{\mathbb{F}}_G , \widehat{\mathbb{F}}_G,\partial_{\c{F}}) \in \b{Op}(\b{sSet}^G)$, which is an operadic variant of Elemendorf's construction of universal spaces \cite[\S 2]{Elm}. The $0$-simplices in $\s{B}_{\c{I}}$ form a discrete, free $G$-operad that contains all operations specified by $\c{I}$, and the remaining simplices kill all homotopy by the extra degeneracy argument.

\end{document}